\documentclass{amsart} 
\usepackage[all]{xy}
\usepackage{amssymb} 
\usepackage{amsmath} 
\usepackage{amscd} 
\usepackage{latexsym} 
\usepackage{amsthm} 
\usepackage{verbatim} 
\usepackage{graphics} 
\usepackage{epsfig} 
\let\cal\mathcal

\def\Iscr{{\cal I}} 
\def\Jscr{{\cal J}}

\def\Mscr{{\cal M}} 
 
\def\Oscr{{\cal O}}

\def\Tscr{{\cal T}} 
 
\def\Vscr{{\cal V}}

\let\blb\mathbb 
 
\def\CC{{\blb C}}

\def \PP{{\blb P}} 
\def \ZZ{{\blb Z}} 
 
\def \NN{{\blb N}} 
\def \RR{{\blb R}}

\def\coker{\operatorname {coker}}

\def\Ext{\operatorname {Ext}}

\def\gr{\operatorname{gr}}

\def\Hilb{\operatorname {Hilb}} 
 
\def\Hom{\operatorname {Hom}}

\def\K0{{K_{0}(\blb P^{2}_{q})}}  
\def\ker{\operatorname {ker}}

\def\P2q{\operatorname {{\blb P}^{2}_{q}}}

\def\pd{\operatorname {pd}}  
\def\pdim{\operatorname {pd}}

\def\PP{\operatorname {\blb P}}

\def\r{\rightarrow}

\def\Tor{\operatorname {Tor}} 

\DeclareMathOperator{\Proj}{Proj}

\newtheorem{lemma}{Lemma}[section] 
\newtheorem{proposition}[lemma]{Proposition} 
\newtheorem{theorem}[lemma]{Theorem} 
 
\newtheorem{condition}[lemma]{Condition}

\newtheorem{lemmas}{Lemma}[subsection] 
\newtheorem{propositions}[lemmas]{Proposition} 
 
\newtheorem{corollarys}[lemmas]{Corollary}

\theoremstyle{definition} 

\newtheorem{example}[lemma]{Example}

{

\newtheorem{case}{Case}

} 

\theoremstyle{remark} 

\newtheorem{remark}[lemma]{Remark} 
\newtheorem{remarks}[lemmas]{Remark}

\newdimen\uboxsep \uboxsep=1ex 
\def\uboxn#1{\vtop to 0pt{\hrule height 0pt depth 0pt\vskip\uboxsep 
\hbox to 0pt{\hss #1\hss}\vss}} 

\def\uboxs#1{\vbox to 0pt{\vss\hbox to 0pt{\hss #1\hss} 
\vskip\uboxsep\hrule height 0pt depth 0pt}}

\numberwithin{equation}{section} 

 

\author{Koen De Naeghel and Michel Van den Bergh}
 
\address{Departement WNI \\
Limburgs Universitair Centrum \\ 
Universitaire Campus \\ 
Building D \\ 
3590 Diepenbeek \\  
Belgium} 

\thanks{The second author is a director of research at the FWO} 
\email[K. De Naeghel]{koen.denaeghel@luc.ac.be} 
\email[M. Van den Bergh]{michel.vandenbergh@luc.ac.be} 
\date{March 17, 2005} 
\title[On incidence between strata of the Hilbert scheme of points on $\PP^{2}$]{On incidence between strata of the Hilbert scheme of points on $\PP^{2}$}

\begin{document} 
\begin{abstract}
The Hilbert scheme of $n$ points in the projective plane has a natural stratification obtained from the associated Hilbert series. In general, the precise inclusion relation between the closures of the strata is still unknown.

In \cite{Guerimand} Guerimand studied this problem for strata whose Hilbert series are as close as possible. Preimposing
a certain technical condition he obtained necessary and sufficient conditions for the incidence of such strata. 

In this paper we present a new approach, based on deformation theory, to Guerimand's result. This allows us to show that the technical condition is not necessary.
\end{abstract}
\maketitle
\tableofcontents

\section{Introduction and main result} 
\label{ref-1-0}
Below $k$ is an algebraically closed field of characteristic zero and $A = k[x,y,z]$. We will consider the Hilbert scheme $\Hilb_{n}(\PP^{2})$ parametrizing zero-dimensional subschemes of length $n$ in $\PP^2$. It is well known that this is a smooth connected projective variety of dimension $2n$. 

Associated to $X \in \Hilb_{n}(\PP^{2})$ there is an ideal $\Iscr_X \subset \Oscr_{\PP^2}$ and a graded ideal $I_X=\oplus_n H^0(\PP^2,\Iscr_X(n))\subset A$. The Hilbert function $h_{X}$ of $X$ is the Hilbert function of the graded ring $A(X)=A/I_X$. Classically $h_X(m)$ is the number of conditions  for a curve of degree $m$ to contain $X$.
Clearly $h_X(m)=n$ for $m\gg 0$.

It seems Castelnuovo was the first to recognize the utility of the difference function (see \cite{Davis})
\[ 
s_X(m) = h_{X}(m) - h_{X}(m-1) 
\] 
Thus $s_X(m)=0$ for $m\gg 0$. Knowing $s_X$ we can reconstruct $h_X$.

It is known \cite{Davis, GMR, GP} that a function $h$ is of the form $h_X$ for $X\in \Hilb_n(\PP^2)$ if and only if $h(m)=0$ for $m<0$ and $h(m)-h(m-1)$ is a so-called \emph{Castelnuovo function of weight $n$}.

A Castelnuovo function \cite{Davis} by definition has the form 
\begin{equation} 
\label{ref-1.1-1} 
s(0)=1,s(1)=2,\ldots,s(\sigma-1)=\sigma \mbox{ and } s(\sigma-1)\ge 
s(\sigma)\ge s(\sigma+1)\ge \cdots \ge 0. 
\end{equation} 
for some integer $\sigma \geq 0$, and the weight of $s$ is the sum of its values.

It is convenient to visualize $s$ using the graph of the staircase function
\[ 
F_{s}: \RR \r \NN: x \mapsto s({\lfloor x \rfloor}) 
\] 
and to divide the area under this graph in unit squares. We will call the result a \emph{Castelnuovo diagram} which, if no confusion arises, we also refer to as $s_{X}$. 

In the sequel we identify a function $f:\ZZ\r \CC$ with its generating function $f(t)=\sum_n f(n) t^n$. We refer to $f(t)$ as a polynomial or a series depending on whether the support of $f$ is finite or not. 

\begin{example} 
$s(t) = 1 + 2t + 3t^{2} + 4t^{3} + 5t^{4} + 5t^{5} + 3t^{6} + 2t^{7} + t^{8} + t^{9} + t^{10}$ is a Castelnuovo polynomial of weight $28$. The corresponding diagram is 
\vspace{0.5cm}

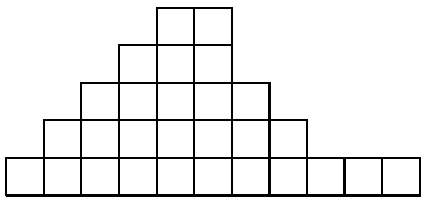
\unitlength 1mm
\begin{picture}(90.00,30.00)(0,0)

\linethickness{0.15mm}
\put(35.00,5.00){\line(1,0){55.00}}

\linethickness{0.15mm}
\put(35.00,10.00){\line(1,0){55.00}}

\linethickness{0.15mm}
\put(40.00,15.00){\line(1,0){35.00}}

\linethickness{0.15mm}
\put(35.00,5.00){\line(0,1){5.00}}

\linethickness{0.15mm}
\put(40.00,5.00){\line(0,1){10.00}}

\linethickness{0.15mm}
\put(45.00,5.00){\line(0,1){15.00}}

\linethickness{0.15mm}
\put(50.00,5.00){\line(0,1){20.00}}

\linethickness{0.15mm}
\put(55.00,5.00){\line(0,1){25.00}}

\linethickness{0.15mm}
\put(90.00,5.00){\line(0,1){5.00}}

\linethickness{0.15mm}
\put(85.00,5.00){\line(0,1){5.00}}

\linethickness{0.15mm}
\put(80.00,5.00){\line(0,1){5.00}}

\linethickness{0.15mm}
\put(75.00,5.00){\line(0,1){10.00}}

\linethickness{0.15mm}
\put(70.00,5.00){\line(0,1){15.00}}

\linethickness{0.15mm}
\put(65.00,5.00){\line(0,1){20.00}}

\linethickness{0.15mm}
\put(45.00,20.00){\line(1,0){25.00}}

\linethickness{0.15mm}
\put(50.00,25.00){\line(1,0){15.00}}

\linethickness{0.15mm}
\put(55.00,30.00){\line(1,0){10.00}}

\linethickness{0.15mm}
\put(60.00,5.00){\line(0,1){25.00}}

\linethickness{0.15mm}
\put(65.00,25.00){\line(0,1){5.00}}

\end{picture}
\end{example} 
We refer to a series $\varphi$ for which $\varphi = h_{X}$ for some $X \in \Hilb_{n}(\PP^{2})$ as a {\em Hilbert function of degree} $n$. The set of all Hilbert functions of degree $n$ (or equivalently the set of all Castelnuovo diagrams of weight $n$) will be denoted by $\Gamma_{n}$.

For $\varphi, \psi \in \Gamma_{n}$ we have that $\psi(t) - \varphi(t)$ is a polynomial, and we write $\varphi \leq \psi$ if its coefficients are non-negative. In this way $\leq$ becomes a partial ordering on $\Gamma_{n}$ and we call the associated directed graph the {\em Hilbert graph}, also denoted by $\Gamma_{n}$.

If $s,t\in \Gamma_n$ are Castelnuovo diagrams such that $s\le t$ then it is
easy to see that $t$ is obtained from $s$ by making  a number of squares ``jump to the
left'' while, at each step, preserving the Castelnuovo property.

\begin{example} 
\label{ref-1.2-2} There are two Castelnuovo diagrams of weight 3.
\vspace{0.5cm}

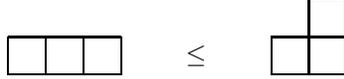
\unitlength 1mm
\begin{picture}(80.00,17.50)(0,0)

\linethickness{0.15mm}
\put(35.00,5.00){\line(0,1){5.00}}

\linethickness{0.15mm}
\put(35.00,10.00){\line(1,0){15.00}}

\linethickness{0.15mm}
\put(50.00,5.00){\line(0,1){5.00}}

\linethickness{0.15mm}
\put(35.00,5.00){\line(1,0){15.00}}

\linethickness{0.15mm}
\put(40.00,5.00){\line(0,1){5.00}}

\linethickness{0.15mm}
\put(45.00,5.00){\line(0,1){5.00}}

\linethickness{0.15mm}
\put(70.00,5.00){\line(0,1){5.00}}

\linethickness{0.15mm}
\put(70.00,5.00){\line(1,0){10.00}}

\linethickness{0.15mm}
\put(80.00,5.00){\line(0,1){10.00}}

\linethickness{0.15mm}
\put(75.00,15.00){\line(1,0){5.00}}

\linethickness{0.15mm}
\put(75.00,5.00){\line(0,1){10.00}}

\linethickness{0.15mm}
\put(70.00,10.00){\line(1,0){10.00}}

\put(60.00,7.50){\makebox(0,0)[cc]{$\leq$}}

\put(30.00,17.50){\makebox(0,0)[cc]{}}

\end{picture}

These distinguish whether three points are collinear or not. The corresponding Hilbert functions are $1,2,3,3,3,3,\ldots \text{ and } 1,3,3,3,3,3,\ldots $. 
\end{example}
\begin{remark}
\label{ref-1.3-3} The number of Castelnuovo diagrams with weight $n$ is equal to the number of partitions of $n$ with distinct parts (or equivalently the number of partitions of $n$ with odd parts) \cite{DV2}. In loc. cit. there is a table of Castelnuovo diagrams of weight up to $6$ as well as some associated data. The Hilbert graph is rather trivial for low values of $n$. The case $n=17$ is more typical (see  Appendix \ref{ref-A-59}).
\end{remark} 
Hilbert functions provide a natural stratification of the Hilbert scheme. For any Hilbert function $\psi$ of degree
$n$ one defines a smooth connected subscheme \cite{DV2,Gotzmann} $H_{\psi}$ of $\Hilb_{n}(\PP^{2})$ by
\[ 
H_{\psi} = \{ X \in \Hilb_{n}(\PP^2) \mid h_{X} = \psi \}. 
\] 
The family $\{ H_{\psi} \}_{\psi \in \Gamma_{n}}$ forms a stratification of $\Hilb_{n}(\PP^{2})$ in the sense that
\[
\overline{H_\psi}\subset \bigcup_{\varphi\le \psi} H_{\varphi}
\]
It follows that if $H_\varphi\subset \overline{H_\psi}$ then $\varphi\le \psi$. The converse implication is in general false and it is still an open problem to find necessary and sufficient conditions for the existence of an inclusion  $H_\varphi\subset \overline{H_\psi}$ \cite{BH, maC, CW, HRW}. This problem is sometimes referred to as the {\em incidence problem}.

Guerimand in his PhD-thesis \cite{Guerimand} introduced two additional necessary conditions for incidence of strata which we now discuss.
\begin{equation} \label{ref-1.2-4}
\mbox{the \emph{dimension condition}: } \dim H_{\varphi} < \dim H_{\psi}  
\end{equation}
This criterion can be used effectively since there are formulas for $\dim H_\psi$ \cite{DV2,Gotzmann}.

The {\em tangent function} $t_{\varphi}$ of a Hilbert function $\varphi \in \Gamma_{n}$ is defined as the Hilbert function of $\Iscr_{X} \otimes_{\PP^2}\Tscr_{\PP^2}$, where $X\in H_\varphi$ is generic. Semi-continuity yields:
\begin{equation} \label{ref-1.3-5} 
\mbox{the \emph{tangent condition}: } t_{\varphi}\ge t_{\psi} 
\end{equation}
Again it is possible to compute $t_\psi$ from $\psi$ (see \cite[Lemme 2.2.4]{Guerimand} and also Proposition \ref{ref-3.3.1-23} below).

\medskip

Let us say that a pair of Hilbert functions $(\varphi, \psi)$ of degree $n$ has {\em length zero} if $\varphi< \psi$ and there are no Hilbert functions $\tau$ of degree $n$ such that $\varphi < \tau < \psi$.\footnote{This is a minor deviation of Guerimand's definition.} It is easy to see $(\varphi, \psi)$ has length zero if and only if the
Castelnuovo diagram of $\psi$ can be obtained from that of $\varphi$ by making a minimal movevement to the left of one square \cite[ Proposition 2.1.7]{Guerimand}.
\begin{example}
Although in the following pair $s_\psi$ is obtained from $s_\varphi$ by moving one square, it is not length zero since $s_{\varphi}$ may be obtained from $s_{\psi}$ by first doing movement $1$ and then $2$. \\
\vspace{-0.8cm}
\unitlength 1mm
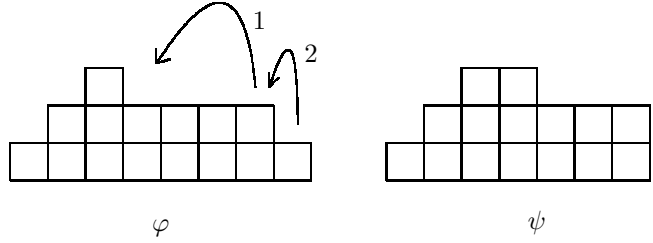
\begin{picture}(105.00,43.13)(0,0)

\linethickness{0.15mm}
\put(20.00,10.00){\line(0,1){5.00}}

\linethickness{0.15mm}
\put(20.00,15.00){\line(1,0){5.00}}

\linethickness{0.15mm}
\put(25.00,15.00){\line(0,1){5.00}}

\linethickness{0.15mm}
\put(25.00,20.00){\line(1,0){5.00}}

\linethickness{0.15mm}
\put(30.00,20.00){\line(0,1){5.00}}

\linethickness{0.15mm}
\put(30.00,25.00){\line(1,0){5.00}}

\linethickness{0.15mm}
\put(40.00,20.00){\line(1,0){15.00}}

\linethickness{0.15mm}
\put(55.00,15.00){\line(0,1){5.00}}

\linethickness{0.15mm}
\put(55.00,15.00){\line(1,0){5.00}}

\linethickness{0.15mm}
\put(60.00,10.00){\line(0,1){5.00}}

\linethickness{0.15mm}
\put(20.00,10.00){\line(1,0){40.00}}

\linethickness{0.15mm}
\put(25.00,10.00){\line(0,1){5.00}}

\linethickness{0.15mm}
\put(30.00,10.00){\line(0,1){10.00}}

\linethickness{0.15mm}
\put(35.00,10.00){\line(0,1){15.00}}

\linethickness{0.15mm}
\put(40.00,10.00){\line(0,1){10.00}}

\linethickness{0.15mm}
\put(45.00,10.00){\line(0,1){10.00}}

\linethickness{0.15mm}
\put(50.00,10.00){\line(0,1){10.00}}

\linethickness{0.15mm}
\put(55.00,10.00){\line(0,1){5.00}}

\linethickness{0.15mm}
\put(25.00,15.00){\line(1,0){30.00}}

\linethickness{0.15mm}
\put(30.00,20.00){\line(1,0){10.00}}

\linethickness{0.15mm}
\put(70.00,10.00){\line(0,1){5.00}}

\linethickness{0.15mm}
\put(70.00,15.00){\line(1,0){5.00}}

\linethickness{0.15mm}
\put(75.00,15.00){\line(0,1){5.00}}

\linethickness{0.15mm}
\put(75.00,20.00){\line(1,0){5.00}}

\linethickness{0.15mm}
\put(80.00,20.00){\line(0,1){5.00}}

\linethickness{0.15mm}
\put(95.00,20.00){\line(1,0){5.00}}

\linethickness{0.15mm}
\put(100.00,15.00){\line(0,1){5.00}}

\linethickness{0.15mm}
\put(75.00,10.00){\line(0,1){5.00}}

\linethickness{0.15mm}
\put(80.00,10.00){\line(0,1){10.00}}

\linethickness{0.15mm}
\put(85.00,10.00){\line(0,1){15.00}}

\linethickness{0.15mm}
\put(90.00,10.00){\line(0,1){15.00}}

\linethickness{0.15mm}
\put(95.00,10.00){\line(0,1){10.00}}

\linethickness{0.15mm}
\put(100.00,10.00){\line(0,1){5.00}}

\linethickness{0.15mm}
\put(105.00,10.00){\line(0,1){5.00}}

\linethickness{0.15mm}
\put(75.00,15.00){\line(1,0){25.00}}

\linethickness{0.15mm}
\put(80.00,20.00){\line(1,0){15.00}}

\put(40.00,3.75){\makebox(0,0)[cc]{$\varphi$}}

\put(90.00,4.38){\makebox(0,0)[cc]{$\psi$}}

\put(40.00,5.00){\makebox(0,0)[cc]{}}

\linethickness{0.15mm}
\put(80.00,25.00){\line(1,0){10.00}}

\linethickness{0.15mm}
\qbezier(52.50,22.50)(50.63,43.13)(39.38,25.63)

\put(53.13,31.25){\makebox(0,0)[cc]{1}}

\linethickness{0.15mm}
\qbezier(58.13,17.50)(58.13,34.38)(54.38,22.50)

\put(60.00,26.88){\makebox(0,0)[cc]{2}}

\linethickness{0.15mm}
\put(54.38,22.50){\line(0,1){1.25}}

\linethickness{0.15mm}
\multiput(54.38,22.50)(0.13,0.13){10}{\line(1,0){0.13}}

\linethickness{0.15mm}
\put(54.38,23.75){\line(0,1){0.63}}

\linethickness{0.15mm}
\put(39.38,25.63){\line(0,1){1.87}}

\linethickness{0.15mm}
\multiput(39.38,25.63)(0.37,0.12){5}{\line(1,0){0.37}}

\linethickness{0.15mm}
\put(70.00,10.00){\line(1,0){30.00}}

\linethickness{0.15mm}
\put(100.00,10.00){\line(1,0){5.00}}

\linethickness{0.15mm}
\put(100.00,15.00){\line(1,0){5.00}}

\linethickness{0.15mm}
\put(100.00,20.00){\line(1,0){5.00}}

\linethickness{0.15mm}
\put(105.00,15.00){\line(0,1){5.00}}

\end{picture}\\

\end{example}
In general a movement of a square by one column is always length zero. A movement by more than one column is length
zero if and only if it is of the form 
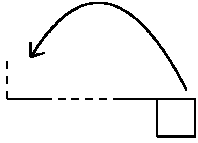
\begin{equation}
\label{ref-1.4-7} 
\unitlength 1mm
\begin{picture}(25.00,31.88)(0,0)

\linethickness{0.15mm}
\put(0.00,10.00){\line(1,0){5.00}}

\linethickness{0.15mm}
\put(15.00,10.00){\line(1,0){10.00}}

\linethickness{0.15mm}
\put(20.00,5.00){\line(0,1){5.00}}

\linethickness{0.15mm}
\put(20.00,5.00){\line(1,0){5.00}}

\linethickness{0.15mm}
\put(25.00,5.00){\line(0,1){5.00}}

\linethickness{0.15mm}
\multiput(0.00,10.00)(0,2.00){3}{\line(0,1){1.00}}

\linethickness{0.15mm}
\multiput(5.00,10.00)(1.82,0){6}{\line(1,0){0.91}}

\linethickness{0.15mm}
\qbezier(23.75,11.25)(13.13,31.88)(3.13,15.63)

\linethickness{0.15mm}
\multiput(3.13,15.63)(0.38,0.13){5}{\line(1,0){0.38}}

\linethickness{0.15mm}
\put(3.13,15.63){\line(0,1){1.88}}

\end{picture}
\end{equation}
The dotted lines represent zero or more squares.

\medskip

The following theorem is the main result of this paper.
\begin{theorem} \label{ref-1.5-8}
Assume that $(\varphi, \psi)$ has length zero. Then $H_{\varphi} \subset \overline{H_{\psi}}$ if and only if 
the dimension condition and the tangent condition hold.
\end{theorem}
This result may be translated into a purely combinatorial (albeit technical) criterion for the existence of an inclusion
$H_\varphi\subset\overline{H_\psi}$ (see Appendix \ref{ref-B-60}).

\medskip

Guerimand proved Theorem \ref{ref-1.5-8} under the additional hypothesis that $(\varphi, \psi)$ is not of ``type zero". A pair of Hilbert series $(\varphi, \psi)$ has {\em type zero} if it is obtained by moving the indicated square in the diagram below.\footnote{It is easy to see that this definition of type zero is equivalent to the one in~\cite{Guerimand}.} \\
\unitlength 1mm
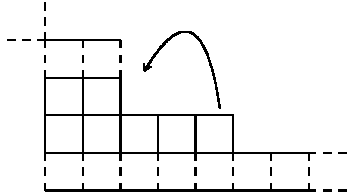
\begin{picture}(90.00,33.44)(0,0)

\linethickness{0.15mm}
\put(50.00,25.00){\line(1,0){10.00}}

\linethickness{0.15mm}
\put(60.00,15.00){\line(0,1){2.50}}

\linethickness{0.15mm}
\put(60.00,15.00){\line(1,0){5.00}}

\linethickness{0.15mm}
\put(65.00,15.00){\line(1,0){10.00}}

\linethickness{0.15mm}
\put(75.00,10.00){\line(0,1){5.00}}

\linethickness{0.15mm}
\put(75.00,10.00){\line(1,0){10.00}}

\linethickness{0.15mm}
\put(60.00,15.00){\line(0,1){5.00}}

\linethickness{0.15mm}
\put(70.00,10.00){\line(0,1){5.00}}

\linethickness{0.15mm}
\put(70.00,10.00){\line(1,0){5.00}}

\linethickness{0.15mm}
\put(70.00,10.00){\line(1,0){5.00}}
\put(70.00,10.00){\line(0,1){5.00}}
\put(75.00,10.00){\line(0,1){5.00}}
\put(70.00,15.00){\line(1,0){5.00}}

\linethickness{0.15mm}
\put(70.00,10.00){\line(1,0){5.00}}
\put(70.00,10.00){\line(0,1){5.00}}
\put(75.00,10.00){\line(0,1){5.00}}
\put(70.00,15.00){\line(1,0){5.00}}

\linethickness{0.15mm}
\qbezier(73.13,15.94)(70.63,33.44)(63.13,20.94)

\put(70.00,30.00){\makebox(0,0)[cc]{}}

\linethickness{0.15mm}
\multiput(45.00,25.00)(2.00,0){3}{\line(1,0){1.00}}

\linethickness{0.15mm}
\multiput(85.00,10.00)(2.00,0){3}{\line(1,0){1.00}}

\linethickness{0.15mm}
\put(63.13,20.94){\line(0,1){0.93}}

\linethickness{0.15mm}
\multiput(63.13,20.94)(0.19,0.13){5}{\line(1,0){0.19}}

\linethickness{0.15mm}
\put(60.00,10.00){\line(0,1){5.00}}

\linethickness{0.15mm}
\put(60.00,10.00){\line(1,0){10.00}}

\linethickness{0.15mm}
\put(65.00,10.00){\line(0,1){5.00}}

\linethickness{0.15mm}
\put(50.00,20.00){\line(1,0){10.00}}

\linethickness{0.15mm}
\put(50.00,10.00){\line(0,1){10.00}}

\linethickness{0.15mm}
\put(50.00,10.00){\line(1,0){10.00}}

\linethickness{0.15mm}
\put(55.00,10.00){\line(0,1){10.00}}

\linethickness{0.15mm}
\put(50.00,15.00){\line(1,0){10.00}}

\linethickness{0.15mm}
\multiput(50.00,20.00)(0,2.00){3}{\line(0,1){1.00}}

\linethickness{0.15mm}
\multiput(55.00,20.00)(0,2.00){3}{\line(0,1){1.00}}

\linethickness{0.15mm}
\multiput(60.00,20.00)(0,2.00){3}{\line(0,1){1.00}}

\linethickness{0.15mm}
\multiput(50.00,5.00)(0,2.00){3}{\line(0,1){1.00}}

\linethickness{0.15mm}
\multiput(55.00,5.00)(0,2.00){3}{\line(0,1){1.00}}

\linethickness{0.15mm}
\multiput(60.00,5.00)(0,2.00){3}{\line(0,1){1.00}}

\linethickness{0.15mm}
\multiput(65.00,5.00)(0,2.00){3}{\line(0,1){1.00}}

\linethickness{0.15mm}
\multiput(70.00,5.00)(0,2.00){3}{\line(0,1){1.00}}

\linethickness{0.15mm}
\multiput(75.00,5.00)(0,2.00){3}{\line(0,1){1.00}}

\linethickness{0.15mm}
\multiput(80.00,5.00)(0,2.00){3}{\line(0,1){1.00}}

\linethickness{0.15mm}
\multiput(85.00,5.00)(0,2.00){3}{\line(0,1){1.00}}

\linethickness{0.15mm}
\multiput(50.00,25.00)(0,2.00){3}{\line(0,1){1.00}}

\linethickness{0.15mm}
\put(50.00,5.00){\line(1,0){35.00}}

\linethickness{0.15mm}
\multiput(85.00,5.00)(2.00,0){3}{\line(1,0){1.00}}

\end{picture}\\
The dotted lines represent zero or more squares.

From the results in  Appendix \ref{ref-B-60} one immediately deduces
\begin{proposition}
Let $\varphi, \psi$ be Hilbert functions of degree $n$ such that $(\varphi, \psi)$ has type zero. Then $H_{\varphi} \subset \overline {H_{\psi}}$.
\end{proposition}
\begin{remark}
The smallest, previously open, incidence problem of type zero seems to be \\

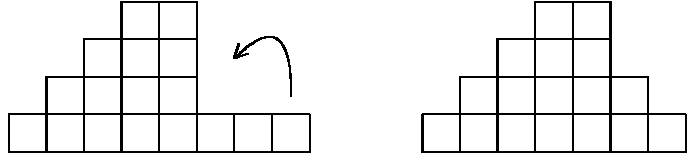
\unitlength 1mm
\begin{picture}(105.00,30.00)(0,0)

\linethickness{0.15mm}
\put(15.00,10.00){\line(0,1){5.00}}

\linethickness{0.15mm}
\put(15.00,15.00){\line(1,0){5.00}}

\linethickness{0.15mm}
\put(20.00,15.00){\line(0,1){5.00}}

\linethickness{0.15mm}
\put(20.00,20.00){\line(1,0){5.00}}

\linethickness{0.15mm}
\put(25.00,20.00){\line(0,1){5.00}}

\linethickness{0.15mm}
\put(25.00,25.00){\line(1,0){5.00}}

\linethickness{0.15mm}
\put(30.00,25.00){\line(0,1){5.00}}

\linethickness{0.15mm}
\put(30.00,30.00){\line(1,0){5.00}}

\linethickness{0.15mm}
\put(35.00,30.00){\line(1,0){5.00}}

\linethickness{0.15mm}
\put(40.00,10.00){\line(0,1){20.00}}

\linethickness{0.15mm}
\put(15.00,10.00){\line(1,0){25.00}}

\linethickness{0.15mm}
\put(40.00,15.00){\line(1,0){15.00}}

\linethickness{0.15mm}
\put(55.00,10.00){\line(0,1){5.00}}

\linethickness{0.15mm}
\put(40.00,10.00){\line(1,0){15.00}}

\linethickness{0.15mm}
\put(45.00,10.00){\line(0,1){5.00}}

\linethickness{0.15mm}
\put(50.00,10.00){\line(0,1){5.00}}

\linethickness{0.15mm}
\put(20.00,15.00){\line(1,0){20.00}}

\linethickness{0.15mm}
\put(25.00,20.00){\line(1,0){15.00}}

\linethickness{0.15mm}
\put(30.00,25.00){\line(1,0){10.00}}

\linethickness{0.15mm}
\put(20.00,10.00){\line(0,1){5.00}}

\linethickness{0.15mm}
\put(25.00,10.00){\line(0,1){10.00}}

\linethickness{0.15mm}
\put(30.00,10.00){\line(0,1){15.00}}

\linethickness{0.15mm}
\put(35.00,10.00){\line(0,1){20.00}}

\linethickness{0.15mm}
\put(70.00,10.00){\line(0,1){5.00}}

\linethickness{0.15mm}
\put(70.00,15.00){\line(1,0){5.00}}

\linethickness{0.15mm}
\put(75.00,15.00){\line(0,1){5.00}}

\linethickness{0.15mm}
\put(75.00,20.00){\line(1,0){5.00}}

\linethickness{0.15mm}
\put(80.00,20.00){\line(0,1){5.00}}

\linethickness{0.15mm}
\put(80.00,25.00){\line(1,0){5.00}}

\linethickness{0.15mm}
\put(85.00,25.00){\line(0,1){5.00}}

\linethickness{0.15mm}
\put(85.00,30.00){\line(1,0){10.00}}

\linethickness{0.15mm}
\put(95.00,10.00){\line(0,1){20.00}}

\linethickness{0.15mm}
\put(70.00,10.00){\line(1,0){25.00}}

\linethickness{0.15mm}
\put(95.00,20.00){\line(1,0){5.00}}

\linethickness{0.15mm}
\put(100.00,10.00){\line(0,1){10.00}}

\linethickness{0.15mm}
\put(100.00,15.00){\line(1,0){5.00}}

\linethickness{0.15mm}
\put(105.00,10.00){\line(0,1){5.00}}

\linethickness{0.15mm}
\put(95.00,10.00){\line(1,0){10.00}}

\linethickness{0.15mm}
\put(75.00,15.00){\line(1,0){25.00}}

\linethickness{0.15mm}
\put(80.00,20.00){\line(1,0){15.00}}

\linethickness{0.15mm}
\put(85.00,25.00){\line(1,0){10.00}}

\linethickness{0.15mm}
\put(90.00,10.00){\line(0,1){20.00}}

\linethickness{0.15mm}
\put(85.00,10.00){\line(0,1){15.00}}

\linethickness{0.15mm}
\put(80.00,10.00){\line(0,1){10.00}}

\linethickness{0.15mm}
\put(75.00,10.00){\line(0,1){5.00}}

\linethickness{0.15mm}
\qbezier(52.50,17.50)(52.50,30.00)(45.00,22.50)

\linethickness{0.15mm}
\multiput(45.00,22.50)(0.38,0.13){5}{\line(1,0){0.38}}

\linethickness{0.15mm}
\multiput(45.00,22.50)(0.13,0.38){5}{\line(0,1){0.38}}

\put(35.00,3.75){\makebox(0,0)[cc]{$\varphi = 1,3,6,10,14,15,16,17,17,\dots$}}

\put(90.00,3.75){\makebox(0,0)[cc]{$\psi = 1,3,6,10,14,16,17,17,\dots$}}

\end{picture}

(see \cite[Exemple A.4.2]{Guerimand}).
\end{remark}
\begin{remark}
Theorem \ref{ref-1.5-8} if false without the condition of $(\varphi,\psi)$ being of length zero. See \cite[Exemple A.2.1]{Guerimand}.
\end{remark}

\bigskip

The authors became interested in the incidence problem while they were studying the deformations of the Hilbert schemes of $\PP^{2}$ which come from non-commutative geometry, see \cite{NS, DV1, DV2}.

It seems that the geometric methods of Guerimand do not apply in a non-commutative context and therefore we developed an alternative approach to the incidence problem based on deformation theory (see \S\ref{ref-2-9}). In this approach the type zero condition turned out to be unnecessary. For this reason we have decided to write down our results first in a purely commutative setting. In a forthcoming paper we will describe the corresponding non-commutative theory.

\section{Outline of the proof of the main theorem}
\label{ref-2-9}
Here and in the rest of this papers we work in the graded category. Thus the notations $\Hom$, $\Ext$ etc\dots never have their ungraded meaning. 

\subsection{Generic Betti numbers}
\label{ref-2.1-10}
Let $X \in \Hilb_{n}(\PP^{2})$. It is easy to see that the graded ideal $I_{X}$ associated to $X$ admits a minimal free resolution of the form
\begin{equation} \label{ref-2.1-11} 
0 \r \oplus_{i}A(-i)^{b_{i}} \r \oplus_{i}A(-i)^{a_{i}} \r I_{X} \r 0 
\end{equation} 
where $(a_{i}),(b_{i})$ are sequences of non-negative integers which have finite support, called the {\em graded Betti numbers} of $I_{X}$ (and $X$). They are related to the Hilbert series of $I_{X}$ as
\begin{equation} \label{ref-2.2-12} 
h_{I_{X}}(t) = h_{A}(t)\sum_{i}(a_{i} - b_{i})t^{i} = \frac{\sum_{i}(a_{i} - b_{i})t^{i}}{(1-t)^{3}}
\end{equation}
So the Betti numbers determine the Hilbert series of $I_{X}$. For generic $X$ (in a stratum $H_\psi$) the converse is true since in that case $a_{i}$ and $b_{i}$ are not both non-zero. We will call such $(a_i)_i$, $(b_i)_i$ \emph{generic Betti numbers}.

\subsection{Four sets of conditions}
We fix a pair of Hilbert series $(\varphi,\psi)$ of length zero. Thus for the associated Castelnuovo functions we have
\begin{equation}
\label{ref-2.3-13}
s_\psi(t)=s_\varphi(t)+t^u-t^{v+1}
\end{equation}
for some integers $0<u\le v$. To prove Theorem \ref{ref-1.5-8} we will show that 4 sets of conditions on $(\varphi,\psi)$ are equivalent.
{\def\thelemma{A}
\begin{condition}
$H_\varphi\subset\overline{H_\psi}$.
\end{condition}
}
{\def\thelemma{B}
\begin{condition}
The dimension and the tangent condition hold for $(\varphi,\psi)$.
\end{condition}
} 
Let $(a_i)_i$ and $(b_i)_i$ be the generic Betti numbers associated to~$\varphi$. The next
technical condition restricts the values of the Betti numbers for $i=u$, $u+1$, $v+2$, $v+3$.
{\def\thelemma{C}
\begin{condition} $a_u\neq 0$, $b_{v+3}\neq 0$ and
\[
\begin{cases}
b_{u+1}\le a_u \le b_{u+1}+1\text{ and } b_{v+3}=a_{v+2}&\\
\text{or} & \text{if $v=u+1$}\\
a_u=b_{u+1}+1 \text{ and } b_{v+3}=a_{v+2}-1&\\
&\\
a_u=b_{u+1}+1\text{ and }b_{v+3}=a_{v+2}&\text{if $v\ge u+2$}
\end{cases}
\]
\end{condition}
}
The last condition is of homological nature. Let $I\subset A$ be a graded ideal corresponding to a generic point of $H_\varphi$. Put
\[
\hat{A}=\begin{pmatrix} A & A \\ 0 & A \end{pmatrix} 
\]
For an ideal $J\subset I$ put 
\[
\hat{J}=\begin{pmatrix} J& I \end{pmatrix}
\]
This is a right $\hat{A}$-module.
{
\def\thelemma{D}
\begin{condition} There exists an ideal $J\subset I$, $h_J(t)=\psi$ such that 
\[
\dim_{k} \Ext^1_{\hat{A}}(\hat{J},\hat{J})< \dim_{k} \Ext^1_{A}(J,J)
\]
\end{condition}
}
In the sequel we will verify the implications
\[
A \Rightarrow B \Rightarrow C \Rightarrow D \Rightarrow A
\]
Here the implication $A\Rightarrow B$ is clear and the implication $B\Rightarrow C$ is purely combinatorial.

\medskip

The implication $C\Rightarrow D$ is based on the observation that $I/J$ must be a so-called truncated point module (see \S\ref{ref-4.1-36} below). This allows us to construct the projective resolution of $J$ from that of $I$ and in this way we can compute $\dim_{k} \Ext^1_A(J,J)$. To compute $\Ext^1_{\hat{A}}(\hat{J},\hat{J})$ we view it as the
tangent space to the moduli-space of pairs $(J,I)$.

\medskip

The implication $D\Rightarrow A$ uses elementary deformation theory. Assume that $D$ holds. Starting from some $\zeta\in \Ext^1_A(J,J)$ (which we view as a first order deformation of $J$), not in the image of $\Ext^1_{\hat{A}}(\hat{J},\hat{J})$ we construct a one-parameter family of ideals  $J_\theta$ such that $J_0=J$ and 
$\pd J_\theta=1$ for $\theta\neq 0$. Since $I$ and $J=J_0$ have the same image in $\Hilb_n(\PP^2)$, this shows that $H_\varphi$ is indeed in the closure of $H_\psi$. 

\section{The implication $B\Rightarrow C$}
\label{ref-3-14}
In this section we translate the length zero condition, the dimension condition and the tangent condition in terms of Betti numbers. As a result we obtain that Condition B implies Condition C.

To make the connection between Betti numbers and Castelnuovo diagrams we frequently use the identities
\begin{equation} \label{ref-3.1-15} 
 \sum_{i \leq l}(a_{i}-b_{i}) =1 + s_{l-1} - s_{l} \quad\mbox{ if }\quad l \geq 0 
\end{equation}
\begin{equation} \label{ref-3.2-16}
a_{l} - b_{l}  = -s_{l}+2s_{l-1}-s_{l-2} \quad \mbox{ if } \quad l > 0
\end{equation}
Throughout we fix a pair of Hilbert functions $(\varphi,\psi)$ of degree $n$ and length zero and we let $s=s_\varphi$,
$\tilde{s}=s_\psi$ be the corresponding Castelnuovo diagrams. Thus we have
\begin{equation}
\label{ref-3.3-17}
\psi(t) = \varphi(t) + t^{u} + t^{u+1} + \cdots + t^{v}
\end{equation}
and
\begin{equation}
\label{ref-3.4-18}
\tilde{s}=s+t^u-t^{v+1}
\end{equation}
for some $0<u\le v$. 

The corresponding generic Betti numbers (cfr \S\ref{ref-2.1-10}) are written as $(a_{i}),(b_{i})$ resp. $(\tilde{a}_{i}),(\tilde{b}_{i})$. We also write
\begin{align*}
\sigma &= \min \{i \mid s_i\ge s_{i+1}\}=\min \{i \mid a_{i} > 0 \}\\
\quad \tilde{\sigma} &= \min \{i \mid \tilde{s}_i\ge
\tilde{s}_{i+1}\}=\ \min \{i \mid \tilde{a}_{i} > 0 \}
\end{align*}

\subsection{Translation of the length zero condition}
The proof of the following result is left to the reader.
\begin{propositions} 
\label{ref-3.1.1-19} If $v\geq u+1$ then we have
\begin{eqnarray*}
\begin{array}{c|ccccccccc}
i & \ldots & u & u+1 & u+2 & \ldots & v+1 & v+2 & v+3 & \ldots \\
\hline 
a_{i} \vrule height 1.2em width 0pt& \ldots & \ast & 0 & 0 & \ldots & 0 & \ast & \ast & \dots \\
b_{i} & \ldots & \ast & \ast & 0 & \ldots & 0 & 0 & \ast &  \ldots 
\end{array}
\end{eqnarray*}
where 
\[
a_{u} \leq b_{u+1}+1, \,\, a_{v+2} > 0, \,\, b_{v+3} \leq a_{v+2}. 
\]
\end{propositions}
This result is based on the identity \eqref{ref-3.2-16}. The zeroes among the Betti numbers are caused by the
 ``plateau'' in $s$ between the $u$'th and the $v+1$'th column (see \eqref{ref-1.4-7}).
 
\subsection{Translation of the dimension condition}
The following result allows us to compare the dimensions of the strata $H_\varphi$ and $H_\psi$.
\begin{propositions} \label{ref-3.2.1-20} 
One has
\begin{equation} 
\dim H_{\psi} = \dim H_{\varphi} + 
\sum_{i=u}^{v}(a_{i}-b_{i}) - \sum_{i=u+3}^{v+3}(a_{i}-b_{i}) + e 
\end{equation} 
and
\begin{equation} \label{ref-3.6-21} 
\begin{split} 
\dim H_{\psi} & = \dim H_{\varphi}-s_{u-2}+s_{u-1}+s_{u+1}-s_{u+2} \\ 
& \hspace{1.5cm}+s_{v-1}-s_{v}-s_{v+2}+s_{v+3}+e \\ 
\end{split} 
\end{equation} 
where 
\begin{eqnarray*} 
e = 
\left\{ 
\begin{array}{cl} 
-1 & \mbox{ if } v = u \\ 
1 & \mbox{ if } v = u+1 \\ 
0 & \mbox{ if } v \geq u+2 
\end{array} 
\right. 
\end{eqnarray*} 
\end{propositions} 
\begin{proof} 
The proof uses only  \eqref{ref-3.4-18}. One has the formula \cite{DV2}
\[ 
\dim H_{\varphi} = 1 + n + c_{\varphi} 
\] 
where $c_{\varphi}$ is the constant term of 
\begin{equation*} 
f_{\varphi}(t) = (t^{-1}-t^{-2})s_{\varphi}(t^{-1})s_{\varphi}(t) 
\end{equation*} 
We find 
\begin{equation*} 
\begin{split} 
f_{\psi}(t) & = (t^{-1}-t^{-2})s_{\psi}(t^{-1})s_{\psi}(t) \\ 
& = (t^{-1}-t^{-2})(s_{\varphi}(t^{-1}) + t^{-u} - t^{-v-1}) 
(s_{\varphi}(t) + t^{u} - t^{v+1}) \\ 
& = (t^{-1}-t^{-2})\biggl(\sum_{i}s_{i}t^{-i} + t^{-u} - t^{-v-1}\biggr) 
\biggl(\sum_{j}s_{j}t^{j} + t^{u} - t^{v+1}\biggr) \\ 
& = f_{\varphi}(t) + (t^{-1}-t^{-2})\biggl(\sum_{i}s_{i}t^{u-i} - 
\sum_{i}s_{i}t^{v+1-i} \\ 
& \hspace{0.7cm} 
+ \sum_{j}s_{j}t^{j-u} - 
\sum_{j}s_{j}t^{j-v-1} - t^{v+1-u} - t^{u-v-1} + 2\biggr) 
\end{split} 
\end{equation*} 
Taking constant terms we obtain \eqref{ref-3.6-21}. Applying \eqref{ref-3.1-15} finishes the proof. 
\end{proof} 
We obtain the following rather strong consequence of the dimension condition.
\begin{corollarys} \label{ref-3.2.2-22}
If $v \geq u+2$ then
\[ 
\dim H_{\varphi} < \dim H_{\psi} \Leftrightarrow a_{u} = b_{u+1}+1 \mbox{ and } a_{v+2} = b_{v+3}
\]
and if this is the case then we have in addition 
\[
\dim H_{\psi} = \dim H_{\varphi} + 1 \mbox{ and } u = \sigma, \,\, a_{u} > 0 \,\, a_{v+2} = b_{v+3} > 0
\]
\end{corollarys} 
\begin{proof} 
Due to Proposition \ref{ref-3.1.1-19} we have $s_{u+1} = s_{u+2}$ and $s_{v-1}=s_{v}$ so \eqref{ref-3.6-21} becomes 
\[ 
\dim H_{\varphi} < \dim H_{\psi} \Leftrightarrow 
(s_{u-2}-s_{u-1})+(s_{v+2}-s_{v+3})<0 
\] 
We have that $1 \leq \sigma \leq u$, which implies $s_{v+2} \geq s_{v+3}$, and either $s_{u-2} \geq s_{u-1}$ or $s_{u-1}=s_{u-2}+1$. From this it is easy to see that we have $(s_{u-2}-s_{u-1})+(s_{v+2}-s_{v+3})<0$ if and only if 
$s_{u-1} = s_{u-2}+1$ and $s_{v+2} = s_{v+3}$.

First assume that this is the case. Then it follows from \eqref{ref-3.1-15} and Proposition \ref{ref-3.1.1-19} that $\sigma = u$ hence $a_{u} > 0$, $b_{u} = 0$. Equation \eqref{ref-3.1-15} together with $s_{u} = s_{u+1}$ gives $\sum_{i \leq u+1}(a_{i} - b_{i}) = 1$ and since $a_{u+1} = 0$ (see Proposition \ref{ref-3.1.1-19}) we have $a_{u} =
b_{u+1} + 1$. Further, \eqref{ref-3.1-15} together with $s_{v+2} = s_{v+3}$ gives $\sum_{i \leq v+3}(a_{i} - b_{i}) = 1$. Combined with $\sum_{i \leq u+1}(a_{i} - b_{i}) = 1$ and Proposition \ref{ref-3.1.1-19} we get $a_{v+2} + (a_{v+3} - b_{v+3}) = 0$ where $a_{v+2} > 0$. This gives $a_{v+2} = b_{v+3} > 0$.

Conversely, assume that $a_{u} = b_{u+1}+1$ and $a_{v+2} = b_{v+3}$. Observe that Proposition \ref{ref-3.1.1-19} implies $s_{u} = s_{u+1}$ and $a_{u+1} = 0$, so using \eqref{ref-3.1-15} yields
\[
1 = \sum_{i \leq u+1}(a_{i} - b_{i}) = \sum_{i \leq u-1}(a_{i} - b_{i}) + a_{u} - b_{u+1}
\]
Since we assumed that $a_{u} = b_{u+1}+1$, we find that $\sum_{i \leq u-1}(a_{i} - b_{i}) = 0$ and using \eqref{ref-3.1-15} again we get $s_{u-2} + 1 = s_{u-1}$.  Next, the fact that $s_{v} = s_{v+1}$ (see Proposition \ref{ref-3.1.1-19}) together with \eqref{ref-3.1-15} yields $\sum_{i \leq v+1}(a_{i} - b_{i}) = 1$.  In combination with
equation \eqref{ref-3.1-15} for $l = v+3$ and Proposition \ref{ref-3.1.1-19} we get that $s_{v+2} - s_{v+3} = a_{v+2} + (a_{v+3} - b_{v+3}) = 0$. Since we assumed that $a_{v+2} = b_{v+3}$ this implies that $s_{v+2} - s_{v+3} = a_{v+3}$. Further, since $b_{v+3} = a_{v+2} > 0$ (see Proposition \ref{ref-3.1.1-19}) we have $a_{v+3} = 0$. We conclude that $s_{u-2} + 1 = s_{u-1}$ and $s_{v+2} = s_{v+3}$ which finishes the proof.
\end{proof} 

\subsection{Translation of the tangent condition}
Recall from the introduction that the tangent function $t_{\varphi}$ is the Hilbert function of $\Iscr_{X} \otimes_{\PP^{2}} \Tscr_{\PP^{2}}$ for $X\in H_\varphi$ generic. 
\begin{propositions} 
\label{ref-3.3.1-23}
(See also \cite[Lemme 2.2.24]{Guerimand}) We have
\begin{equation} \label{ref-3.7-24} 
t_{\varphi}(t) = h_{\Tscr_{\PP^2}}(t) - (3t^{-1}-1)\varphi(t) + 
\sum_{i}b_{i+3}t^{i} 
\end{equation} 
\end{propositions}
\begin{proof} 
From the exact sequence
\[ 
0\r \Tscr_{\PP^{2}} \r \Oscr(2)^{3}\r \Oscr(3) \r 0 
\] 
we deduce
\begin{equation} \label{ref-3.8-25}
H^1(\PP^2,\Tscr_{\PP^2}(n))=
\begin{cases}
k&\text{if $n=-3$}\\
0&\text{otherwise}
\end{cases}
\end{equation}
Let $\Iscr=\Iscr_X$ ($X$ generic) and consider the associated resolution.
\begin{equation*} 
0 \r \oplus_{j} \Oscr(-j)^{b_{j}} \r \oplus_{i} \Oscr(-i)^{a_{i}} \r 
\Iscr \r 0 
\end{equation*} 
Tensoring with $\Tscr_{\PP^2}(n)$ and applying the long exact sequence for $H^\ast(\PP^2,-)$ we obtain an exact sequence 
\begin{multline*}
0 \r \oplus_{j} \Gamma(\PP^2,\Tscr(n-j)^{b_{j}} )\r 
\oplus_{i} \Gamma(\PP^2,\Tscr(n-i)^{a_{i}}) \r 
\Gamma(\PP^2,\Iscr\otimes \Tscr(n))\r \\
\oplus_{j} H^1(\PP^2,\Tscr(n-j)^{b_{j}} )\r 
\oplus_{i} H^1(\PP^2,\Tscr(n-i)^{a_{i}})
\end{multline*}
It follows from \eqref{ref-3.8-25} that the rightmost arrow is zero. This easily yields the required formula. 
\end{proof} 
\begin{remarks} 
The previous proposition has an easy generalization which is perhaps useful and which is proved in the same way. Let $M$ be the second syzygy of a finite dimensional graded $A$-module $F$ and let $\Mscr$ be the associated coherent sheaf. Write $h_M(t)=q_M(t)/(1-t)^3$. Then the Hilbert series of $\Iscr_X\otimes\Mscr$ is given by
\[
q_M(t)h_{I_X}(t)+h_{\Tor_1^A(F,I_X)}(t)
\]
The case where $\Mscr$ is the tangent bundle corresponds to $F=k(3)$.
\end{remarks}
\begin{propositions} \label{ref-3.3.3-26} 
We have
\begin{enumerate} 
\item 
$t_{\psi}(l) \leq t_{\varphi}(l)$ for $l \neq u-3,v$ 
\item 
$t_{\psi} \leq t_{\varphi} \Leftrightarrow a_{u} \neq 0 \mbox{ and } b_{v+3} \neq 0$ 
\end{enumerate} 
\end{propositions} 
\begin{proof} 
The proof uses only \eqref{ref-3.4-18}. Comparing \eqref{ref-3.7-24} for $\varphi$ and $\psi$ gives 
\begin{equation} \label{ref-3.9-27} 
t_{\varphi}(t) - t_{\psi}(t) = 3t^{u-1} + 2(t^{u} + t^{u+1} + \ldots + t^{v-1}) - t^{v} + \sum_{i}(b_{i+3} - \tilde{b}_{i+3})t^{i} 
\end{equation} 
where we have used \eqref{ref-3.3-17}. In order to prove the statements, we have to estimate the polynomial 
$\sum_{i}(b_{i+3} - \tilde{b}_{i+3})t^{i}$. For this, substituting \eqref{ref-2.2-12} for $\varphi$ and $\psi$ in 
\eqref{ref-3.3-17} gives 
\begin{equation*} 
\begin{split} 
\sum_{i}(\tilde{a}_{i} - \tilde{b}_{i})t^{i} & = \sum_{i}(a_{i} - b_{i})t^{i}- (t^{u}-t^{v+1})(1-t)^{2} \\ 
& = \sum_{i}(a_{i} - b_{i})t^{i} - t^{u} + 2t^{u+1} -t^{u+2} + t^{v+1} -2t^{v+2} + t^{v+3} 
\end{split} 
\end{equation*} 
hence 
\begin{equation} \label{ref-3.10-28} 
\begin{split} 
\tilde{a}_{u} - \tilde{b}_{u}  & = a_u - b_u - 1 \\ 
\tilde{a}_{u+1} - \tilde{b}_{u+1} & = a_{u+1} - b_{u+1} + 
\left\{ 
\begin{array}{cl} 
3 & \mbox{ if } v = u \\
2 & \mbox{ if } v \geq u+1  
\end{array} 
\right. 
\\ 
\tilde{a}_{u+2} - \tilde{b}_{u+2} & = a_{u+2} - b_{u+2} + 
\left\{ 
\begin{array}{cl} 
-3 & \mbox{ if } v = u \\ 
0 & \mbox{ if } v = u+1 \\ 
-1 & \mbox{ if } v \geq u+2 
\end{array} 
\right. 
\\ 
\tilde{a}_{v+1} - \tilde{b}_{v+1} & = a_{v+1} - b_{v+1} + 
\left\{ 
\begin{array}{cl} 
3 & \mbox{ if } v = u \\ 
0 & \mbox{ if } v = u+1 \\ 
1 & \mbox{ if } v \geq u+2 
\end{array} 
\right. 
\\ 
\tilde{a}_{v+2} - \tilde{b}_{v+2} & = a_{v+2} - b_{v+2} + 
\left\{ 
\begin{array}{cl} 
-3 & \mbox{ if } v = u \\
-2 & \mbox{ if } v \geq u+1  
\end{array} 
\right. 
\\ 
\tilde{a}_{v+3} - \tilde{b}_{v+3}  & = a_{v+3} - b_{v+3} + 1 \\ 
\tilde{a}_{l} - \tilde{b}_{l}  & = a_l - b_l  \quad \mbox{ if } l \not\in \{ u,u+1,u+2,v+1,v+2,v+3 \} 
\end{split} 
\end{equation} 
To obtain information about the differences $b_{i+3} - \tilde{b}_{i+3}$, we observe that for $c \geq 0$ and for all integers $l$ we have 
\begin{equation} \label{ref-3.11-29} 
\begin{split} 
\tilde{a}_{l} - \tilde{b}_{l}  & = a_l - b_l + c  \Rightarrow \tilde{b}_{l} 
\leq b_{l} \\ 
\tilde{a}_{l} - \tilde{b}_{l}  & = a_l - b_l - c  \Rightarrow \tilde{b}_{l} 
\leq b_{l} + c 
\end{split} 
\end{equation} 
Indeed, first let $\tilde{a}_{l} - \tilde{b}_{l}  = a_l - b_l + c$. In case $0 \leq b_{l} \leq c$ then 
$0 = \tilde{b}_{l} \leq b_{l}$. And in case $c < b_{l}$ then $\tilde{b}_{l} = b_{l}-c \leq b_{l}$. \\ 
Second, let $\tilde{a}_{l} - \tilde{b}_{l}  = a_l - b_l - c$. In case $0 \leq a_{l} \leq c$ then $\tilde{a}_{l} = 0$ hence $\tilde{b}_{l}  = b_{l} + c - a_{l} \leq b_{l} + c$. And in case $c < a_{l}$ then $0 = \tilde{b}_{l} \leq c = b_{l} + c$. So this proves \eqref{ref-3.11-29}. \\ 

Applying \eqref{ref-3.11-29} to \eqref{ref-3.10-28} yields 
\begin{equation} \label{ref-3.12-30} 
\begin{split} 
\tilde{b}_{u} & \leq b_{u} + 1 \\ 
\tilde{b}_{u+1} & \leq b_{u+1} \\ 
\tilde{b}_{u+2} & \leq b_{u+2} + 
\left\{ 
\begin{array}{cl} 
3 & \mbox{ if } v = u \\ 
0 & \mbox{ if } v = u+1 \\ 
1 & \mbox{ if } v \geq u+2 
\end{array} 
\right. 
\\ 
\tilde{b}_{v+1} & \leq b_{v+1} \\ 
\tilde{b}_{v+2} & \leq b_{v+2} + 
\left\{ 
\begin{array}{cl} 
3 & \mbox{ if } v = u \\ 
2 & \mbox{ if } v \geq u+1 
\end{array} 
\right. 
\\ 
\tilde{b}_{v+3} & \leq b_{v+3} \\ 
\tilde{b}_{l} & \leq b_{l} \mbox{ if } l \not\in \{ u,u+1,u+2,v+1,v+2,v+3 \} 
\end{split}  
\end{equation} 
Now we are able to prove the first statement. Combining \eqref{ref-3.12-30} and \eqref{ref-3.9-27} gives 
\begin{equation} 
\begin{split} 
t_{\varphi}(t) - t_{\psi}(t) \geq 
\left \{ 
\begin{array}{ll} 
-t^{u-3} - t^{v} & \mbox{ if } v=u \\ 
-t^{u-3} + 3t^{u-1} - t^{v} & \mbox{ if } v=u+1 \\ 
-t^{u-3} + 2(t^{u-1} + t^{u} + \ldots + t^{v-2}) - t^{v} & 
\mbox{ if } v \geq u+2 
\end{array} 
\right. 
\end{split} 
\end{equation} 
and therefore $t_{\varphi}(t) - t_{\psi}(t) \geq -t^{u-3} - t^{v}$ which concludes the proof of the first statement. \\ 

For the second part, assume that $t_{\psi} \leq t_{\varphi}$. Equation \eqref{ref-3.9-27} implies that 
\begin{equation} \label{ref-3.14-31} 
\begin{split} 
\tilde{b}_{u} & \leq b_{u} \\ 
\tilde{b}_{v+3} & \leq b_{v+3} - 1 
\end{split} 
\end{equation} 
Since $\tilde{b}_{v+3} \geq 0$ we clearly have $b_{v+3} > 0$. Assume, by contradiction, that $a_{u} = 0$. From \eqref{ref-3.10-28} we have $\tilde{a}_{u} - \tilde{b}_{u}  = a_u - b_u - 1$ hence $\tilde{a}_{u} = 0$ and $\tilde{b}_{u} = b_u + 1$. But this gives a contradiction with \eqref{ref-3.14-31}. Therefore 
\[ 
t_{\psi} \leq t_{\varphi} \Rightarrow a_{u} > 0 \mbox{ and } 
b_{v+3} > 0 
\] 
To prove the converse let $a_{u} > 0$ and $b_{v+3} > 0$. Due to the first part we only need to prove that $t_{\psi}(u-3) \leq t_{\varphi}(u-3)$ and $t_{\psi}(v) \leq t_{\varphi}(v)$. Equation \eqref{ref-3.9-27} gives us
\begin{equation} \label{ref-3.15-32} 
\begin{split} 
& t_{\varphi}(u-3) - t_{\psi}(u-3) = b_{u} - \tilde{b}_{u} \\ 
& t_{\varphi}(v) - t_{\psi}(v) = b_{v+3} - \tilde{b}_{v+3} - 1 
\end{split} 
\end{equation} 
while from \eqref{ref-3.10-28} we have 
\begin{equation*} 
\begin{split} 
& \tilde{a}_{u} - \tilde{b}_{u}  = a_u - b_u - 1 \\ 
& \tilde{a}_{v+3} - \tilde{b}_{v+3}  = a_{v+3} - b_{v+3} + 1 
\end{split} 
\end{equation*} 
Since $a_{u} > 0$, $b_{v+3} > 0$ we have $b_{u} = 0$, $a_{v+3} = 0$ hence 
\begin{equation*} 
\begin{split} 
& \tilde{a}_{u} - \tilde{b}_{u}  = a_u - 1 \\ 
& \tilde{a}_{v+3} - \tilde{b}_{v+3}  = - b_{v+3} + 1 
\end{split} 
\end{equation*} 
which implies $\tilde{a}_{u} - \tilde{b}_{u} \geq 0$, $\tilde{a}_{v+3} - \tilde {b}_{v+3} \leq 0$ hence $\tilde{b}_{u} = 0$, $\tilde{a}_{v+3} = 0$. Thus $b_{u} = \tilde{b}_{u} = 0$ and $\tilde {b}_{v+3}  =  b_{v+3} - 1$. Combining with \eqref{ref-3.15-32} this proves that $t_{\varphi}(u-3) = t_{\psi}(u-3)$ and $t_{\varphi}(v) = t_{\psi}(v)$, finishing the proof. 
\end{proof} 

\subsection{Combining everything}
In this section we prove that Condition $B$ implies Condition $C$. So assume that Condition $B$ holds.

Since the tangent condition holds we have by Proposition \ref{ref-3.3.3-26}
\[
a_u\neq 0\qquad\text{and}\qquad b_{v+3}\neq 0
\]
This means there is nothing to prove if $u=v$. We discuss the two remaining cases.

\begin{case}
$v = u+1$
\end{case}
The fact that $a_{u} \neq 0$, $b_{v+3} \neq 0$ implies $b_{u} = 0$, $a_{v+3} = 0$. Proposition \ref{ref-3.2.1-20} combined with Proposition \ref{ref-3.1.1-19} now gives
\[
\dim H_{\psi} = \dim H_{\varphi} + a_{u} - b_{u+1} - a_{v+2} + b_{v+3}
+1
\]
Hence $0 \leq (a_{u} - b_{u+1}) + (b_{v+3} - a_{v+2})$. But Proposition \ref{ref-3.1.1-19} also states that $a_{u} \leq b_{u+1}+1$, $a_{v+2} > 0$ and $b_{v+3} \leq a_{v+2}$. Therefore either we have that
\[ 
b_{u+1} \leq a_{u} \leq b_{u+1}+1 \mbox{ and } b_{v+3} = a_{v+2}
\]
or
\[
a_{u} = b_{u+1}+1 \mbox{ and } b_{v+3} = a_{v+2}-1
\]
\begin{case}
$v \geq u+2$
\end{case}
It follows from Corollary \ref{ref-3.2.2-22} that
\[
a_{u} = b_{u+1} + 1 \mbox{ and } a_{v+2} = b_{v+3}
\]  
This finishes the proof.
\begin{remarks} 
\label{ref-3.4.1-33} The reader will have noticed that our proof of the implication $B\Rightarrow C$ is rather involved.
Since the equivalence of $B$ and $C$ is purely combinatorial it can be checked directly for individual $n$. Using a computer we have verified the equivalence of $B$ and $C$ for $n\le 70$.

As another independent verification we have a direct proof of the implication $C\Rightarrow B$ (i.e.\ without going through the other conditions).
\end{remarks}
\begin{remarks}
The reader may observe that in case $v = u$ we have
\begin{eqnarray}
\label{ref-3.16-34}
t_{\psi} \leq t_{\varphi} \Rightarrow \dim H_{\varphi} < \dim H_{\psi}
\end{eqnarray}
while if $v \geq u+2$ we have
\begin{eqnarray}
\label{ref-3.17-35}
\dim H_{\varphi} < \dim H_{\psi} \Rightarrow t_{\psi} \leq t_{\varphi}
\end{eqnarray}
It is easy to construct counter examples which show that the reverse implications do not hold, and neither \eqref{ref-3.16-34} nor \eqref{ref-3.17-35} is valid in case $v = u+1$.
\end{remarks}
\section{The implication $C\Rightarrow D$}
In this section $(\varphi,\psi)$ will have the same meaning as in \S\ref{ref-3-14} and we also keep the associated notations. 

\subsection{Truncated point modules}
\label{ref-4.1-36}
A truncated point module of length $m$ is a graded $A$-module generated in degree zero with Hilbert series $1 + t + \cdots + t^{m-1}$.

If $F$ is a truncated point module of length $>1$ then there are two independent homogeneous linear forms $l_{1},l_{2}$
vanishing on $F$ and their intersection defines a point $p\in \PP^2$. We may choose basis vectors $e_i\in F_i$ such that 
\[
xe_i=x_pe_{i+1}, \qquad ye_i=y_pe_{i+1}, \quad ze_i=z_pe_{i+1}
\]
where $(x_p,y_p,z_p)$ is a set of homogeneous coordinates of $p$. It follows that if $f\in A$ is homogeneous of degree $d$ and $i+d\le m-1$ then
\[
fe_i = f_p e_{i+d}
\]
where $(-)_p$ stands for evaluation in  $p$ (with respect to the homogeneous coordinates $(x_p,y_p,z_p)$).

If $G=\oplus_i A(-i)^{c_i}$ then we have
\begin{equation}
\label{ref-4.1-37}
\Hom_A(G,F)=\oplus_{0\le i \le m-1} F_i^{c_i}\cong k^{\sum_{0\le i\le m-1} c_i}
\end{equation}
where the last identification is made using the basis $(e_i)_i$ introduced above. 

\medskip

In the sequel we will need the minimal projective resolution of a truncated point module $F$ of length $m$. It is easy to see that it is given by 
{\tiny
\begin{equation}
\label{ref-4.2-38}
0\r A(-m-2)
\xrightarrow[f_3]{
\begin{pmatrix}
l_1\\
l_2\\
\rho
\end{pmatrix}\cdot
}
A(-m-1)^2\oplus A(-2)
\xrightarrow[f_2]{
\begin{pmatrix}
0& -\rho & l_2\\
\rho& 0 & -l_1\\
-l_2 & l_1 & 0
\end{pmatrix}\cdot
}
A(-1)^2\oplus A(-m)
\xrightarrow[f_1]{\begin{pmatrix}l_1& l_2&\rho\end{pmatrix}\cdot}
A\r F\r 0
\end{equation}
}
where $l_1,l_2$ are the linear forms vanishing on $F$ and $\rho$ is a form of degree $m$ such that $\rho_p\neq 0$ for the point $p$ corresponding to $F$. Without loss of generality we may and we will assume that $\rho_p=1$.
 
\subsection{A complex whose homology is  $J$}
\label{ref-4.2-39}
In this section $I$ is a graded ideal corresponding to a generic point in $H_\varphi$. The following lemma gives the connection between truncated point modules and Condition D.
\begin{lemmas} 
If an ideal $J\subset I$ has Hilbert series $\psi$ then $I/J$ is a (shifted by grading) truncated point module of length $v+1-u$.
\end{lemmas}
\begin{proof} 
Since $F=I/J$ has the correct (shifted) Hilbert function, it is sufficient to show that $F$ is generated in degree $u$.

If $v=u$ then there is nothing to prove. If $v\ge u+1$ then by Proposition \ref{ref-3.1.1-19} the generators of $I$ are in degrees $\le u$ and $\ge v+2$. Since $F$ lives in degrees $u,\ldots,v$ this proves what we want.
\end{proof}
Let $J,F$ be as in the previous lemma. Below we will need a complex whose homology is $J$. We write the minimal resolution of $F$ as
\[
0\r G_3 \xrightarrow{f_3} G_2 \xrightarrow{f_2} G_1 \xrightarrow{f_1} G_0 \xrightarrow{} F \r  0
\]
where the maps $f_{i}$ are as in \eqref{ref-4.2-38}, and the minimal resolution of $I$ as
\[
0\r F_1 \r F_0 \r I \r 0 
\]
The map $I\r F$ induces a map of projective resolutions {\small
\begin{eqnarray} \label{ref-4.3-40} 
\begin{CD} 
&& && 0 @>>> F_1 @>{M}>> F_0 @>>> I @>>> 0 \\ 
&& && && @V{\gamma_1}VV @V{\gamma_0}VV @VVV && \\ 
0 @>>> G_3 @>{f_{3}}>> G_2 @>{f_2}>> G_1 @>{f_1}>> G_0 @>{f_0}>> F @>>> 0 
\end{CD} 
\end{eqnarray} 
}  
Taking cones yields that $J$ is the  homology at $G_1\oplus F_0$ of the following complex
\begin{equation}
\label{ref-4.4-41}
0 \r G_3 \xrightarrow{
\begin{pmatrix}f_3\\ 0\end{pmatrix}
}G_2\oplus F_1
\xrightarrow{
\begin{pmatrix}
f_2 & \gamma_1\\
0 & -M
\end{pmatrix}
}
G_1\oplus F_0
\xrightarrow{\begin{pmatrix}f_1 & \gamma_0\end{pmatrix}}
G_0\r 0
\end{equation}
Note that the rightmost map is split here. By selecting an explicit splitting we may construct a free resolution of $J$, but it will be convenient not to do this.

\medskip

For use below we note that the map $J\r I$ is obtained from taking homology of the following map of complexes.
\begin{equation}
\label{ref-4.5-42}
\xymatrix{
0 \ar[r] & G_3 \ar[r]^(0.4)
{\begin{pmatrix}f_3\\ 0\end{pmatrix}}
 & G_2\oplus F_1 \ar[d]_{\begin{pmatrix} 0 & -1\end{pmatrix}}
\ar[rr]^{
\begin{pmatrix}
f_2 & \gamma_1\\
0 & -M
\end{pmatrix}
}
&&
G_1\oplus F_0\ar[d]^{\begin{pmatrix} 0 & 1\end{pmatrix}}
\ar[rr]^(0.6){\begin{pmatrix}f_1 & \gamma_0\end{pmatrix}}
&&G_0 \ar[r]& 0\\
&0\ar[r]& F_1 \ar[rr]_M &&F_0\ar[r]& 0
}
\end{equation}

\subsection{The Hilbert scheme of an ideal}
In this section $I$ is a graded ideal corresponding to a generic point in $H_\varphi$. 

Let $\Vscr$ be the Hilbert scheme of graded quotients $F$ of $I$ with Hilbert series $t^u+\cdots+t^v$. To see that $\Vscr$ exists one may realize it as a closed subscheme of 
\[
\Proj S(I_u\oplus\cdots \oplus I_v)
\]
where $SV$ is the symmetric algebra of a vector space $V$. Alternatively see \cite{AZ2}.

\medskip

We will give an explicit description of $\Vscr$ by equations. Here and below we use the following convention: if $N$ is a matrix with coefficients in $A$ representing a map $\oplus_j A(-j)^{d_j}\r \oplus_i A(-i)^{c_i}$ then $N(p,q)$ stands for the submatrix of $N$ representing the induced map $ A(-q)^{d_q}\r  A(-p)^{c_p}$.

\medskip

We now distinguish two cases.

\begin{itemize}
\item[$v=u$]  
In this case it is clear that $\Vscr\cong \PP^{a_u-1}$.
\item[$v\ge u+1$] 
Let $F\in \Vscr$ and let $p\in \PP^2$ be the associated point. Let $(e_i)_i$ be a basis for $F$ as in \S\ref{ref-4.1-36}. The map $I\r F$ defines a map 
\[
\lambda:A(-u)^{a_u}\r F
\]
such that the composition
\begin{equation}
\label{ref-4.6-43}
A(-u-1)^{b_{u+1}} \xrightarrow{M(u,u+1)\cdot} A(-u)^{a_u} \r F
\end{equation}
is zero. 

We may view $\lambda$ as a scalar row vector as in \eqref{ref-4.1-37}. The fact that \eqref{ref-4.6-43} has zero composition then translates into the  condition
\begin{equation}
\label{ref-4.7-44}
\lambda \cdot M(u,u+1)_p=0
\end{equation}
It is easy to see that this procedure is reversible and that the equations \eqref{ref-4.7-44} define $\Vscr$ as a subscheme of $\PP^{a_u-1}\times \PP^2$.
\end{itemize}
\begin{propositions} 
\label{ref-4.3.1-45} Assume that Condition C holds. Then $\Vscr$ is smooth and
\[
\dim \Vscr=
\begin{cases}
a_u-1 &\text{if $v=u$} \\
a_u+1-b_{u+1} &\text{if $v\ge u+1$}
\end{cases}
\]
\end{propositions}
\begin{proof}
The case $v=u$ is clear so assume $v\ge u+1$. If we look carefully at \eqref{ref-4.7-44} then we see that it describes
$\Vscr$ as the zeroes of $b_{u+1}$ generic sections in the very ample line bundle $\Oscr_{\PP^{a_u-1}}(1) \boxtimes \Oscr_{\PP^2}(1)$ on $\PP^{a_u-1}\times \PP^2$. It follows from Condition C that $b_{u+1}\le \dim (\PP^2\times \PP^{a_u-1})=a_u+1$. Hence by Bertini (see \cite{H}) we deduce that $\Vscr$ is smooth of dimension $a_u+1-b_{u+1}$.
\end{proof}
\subsection{Estimating the dimension of {$\Ext^1_A(J,J)$}}
In this section $I$ is a graded ideal corresponding to a generic point of $H_\varphi$. We prove the following result
\begin{propositions} 
\label{ref-4.4.1-46}
Assume that Condition C holds. Then there exists $F\in \Vscr$ such that for $J=\ker(I\r F)$ we have
\begin{equation}
\label{ref-4.8-47}
\dim_{k} \Ext^1_A(J,J)
\ge
\begin{cases}
\dim H_\psi + a_{v+3} = \dim H_\psi & \text{if $v=u$} \\
\dim H_\psi + a_{v+2} - b_{v+3} + 1 & \text{if $v=u+1$} \\
\dim H_\psi + a_{v+2} - b_{v+3} + 2=\dim H_\psi + 2 & \text{if $v \ge u+2$}
\end{cases}
\end{equation}
\end{propositions}
It will become clear from the proof below that in case $v\ge u+1$ the righthand side of \eqref{ref-4.8-47} is one higher than the expected dimension.

\medskip

Below let $J \subset I$ be an arbitrary ideal such that $h_J = \psi$. Put $F = I/J$.
\begin{propositions} 
We have
\[
\dim_{k} \Ext^1_A(J,J) = \dim H_\psi + \dim_{k} \Hom_A(J,F(-3))
\]
\end{propositions}
\begin{proof}
For $M,N\in \gr A$ write
\[
\chi(M,N)=\sum_i (-1)^i \dim_{k} \Ext^i_A(M,N)
\]
Clearly $\chi(M,N)$ only depends on the Hilbert series of $M$, $N$. Hence, taking  $J'$ to be an arbitrary point in $H_\psi$ we have
\[
\chi(J,J) = \chi(J',J') = 1 - \dim_{k} \Ext^1_A(J',J') = 1 - \dim H_\psi
\]
where in the third equality we have used that $\Ext^1_A(J',J')$ is the tangent space to $H_\psi$ \cite{DV2}.

Since $J$ has no socle we have $\pdim J \le 2$. Therefore $\Ext^i_A(J,J) = 0$ for $i \ge 3$. It follows that
\begin{align*}
\dim_{k} \Ext^1_A(J,J) & = -\chi(J,J) + 1 + \dim_{k} \Ext^2_A(J,J) \\
& = \dim H_\psi + \dim_{k} \Ext^3_A(F,J)
\end{align*}
By the approriate version of Serre duality we have
\[
\Ext^3_A(F,J) = \Hom_A(J,F \otimes \omega_A)^\ast = \Hom_A(J,F(-3))^\ast
\]
This finishes the proof.
\end{proof}
\begin{proof}[Proof of Proposition \ref{ref-4.4.1-46}] 
It follows from the previous result that we need to control $\dim_{k} \Hom_A(J,F(-3))$. Of course we assume throughout that Condition C holds and we also use Proposition \ref{ref-3.1.1-19}.
\setcounter{case}{0}
\begin{case}
Assume $v=u$. For degree reasons any extension between $F$ and $F(-3)$ must be split. Thus we have $\Hom_A(F,F(-3)) = \Ext^1_A(F,F(-3)) = 0$. Applying $\Hom_A(-,F(-3))$ to 
\[
0 \r J \r I \r F \r 0
\]
we find
\[
\Hom_A(J,F(-3)) = \Hom_A(I,F(-3))
\]
Hence
\[
\dim_{k}\Hom_A(J,F(-3)) = a_{v+3} = 0
\]
\end{case}
\begin{case} 
Assume $v=u+1$. As in the previous case we find $\Hom_A(J,F(-3)) = \Hom_A(I,F(-3))$.

Thus a map $J\r F(-3)$ is now given (using Proposition \ref{ref-3.1.1-19}) by a map
\[
\beta:A(-v-2)^{a_{v+2}}\r F(-3)
\]
(identified with a scalar vector as in \eqref{ref-4.1-37}) such that the composition
\[
A(-v-3)^{b_{v+3}} \xrightarrow{M(v+2,v+3)} A(-v-2)^{a_{v+2}} \xrightarrow{\beta} F(-3)
\]
is zero. This translates into the condition
\begin{equation}
\label{ref-4.9-48}
\beta\cdot M(v+2,v+3)_p=0
\end{equation}
where $p$ is the point corresponding to $F$. Now $M(v+2,v+3)$ is a $a_{v+2}\times b_{v+3}$ matrix. Since $b_{v+3}\le a_{v+2}$ (by Proposition \ref{ref-3.1.1-19}) we would expect \eqref{ref-4.9-48} to have $a_{v+2}-b_{v+3}$ independent solutions. To have more, $M(v+2,v+3)$ has to have non-maximal rank.  I.e.\ there should be a non-zero solution to the equation
\begin{equation}
\label{ref-4.10-49}
M(v+2,v+3)_p\cdot \delta=0
\end{equation}
This should be combined with (see \eqref{ref-4.7-44})
\begin{equation}
\label{ref-4.11-50}
\lambda\cdot M(u,u+1)_p=0
\end{equation}
We view \eqref{ref-4.10-49} and \eqref{ref-4.11-50} as a system of $a_{v+2}+b_{u+1}$ equations in $\PP^{a_{u}-1}\times \PP^2\times \PP^{b_{v+3}-1}$.

Since (Condition C)
\[
a_{v+2}+b_{u+1}\le \dim (\PP^{a_{u}-1}\times \PP^2\times \PP^{b_{v+3}-1})=a_u+b_{v+3}
\]
the system \eqref{ref-4.10-49}\eqref{ref-4.11-50} has a solution provided the
divisors in $\PP^{a_{u}-1}\times \PP^2\times \PP^{b_{v+3}-1}$ determined by the
equations of the system have non-zero intersection product.

Let $r,s,t$ be the hyperplane sections in $\PP^{a_{u}-1}$, $\PP^2$ and
$\PP^{b_{v+3}-1}$ respectively. The Chow ring of $\PP^{a_{u}-1}\times
\PP^2\times \PP^{b_{v+3}-1}$ is given by
\begin{equation}
\label{ref-4.12-51}
\ZZ[r,s,t]/(r^{a_u}, s^3, t^{b_{v+3}})
\end{equation}
The intersection product we have to compute is 
\[
(s+t)^{a_{v+2}}(r+s)^{b_{u+1}}
\]
This product contains the terms
\begin{gather*}
t^{a_{v+2}-2}s^2 r^{b_{u+1}}\\
t^{a_{v+2}-1}s^2 r^{b_{u+1}-1}\\
t^{a_{v+2}}s^2 r^{b_{u+1}-2}
\end{gather*}
at least one of which is non-zero in \eqref{ref-4.12-51} (using Condition C).
\end{case}
\begin{case}
Now assume $v\ge u+2$. We compute $\Hom(J,F(-3))$ as the homology of $\Hom_A((\mathrm{eq. }\ref{ref-4.4-41}),F(-3))$. Since $G_0=A(-u)$ we have $\Hom_A(G_0,F(-3))=0$ and hence a map $J\r F(-3)$ is given by a map
\[
G_1\oplus F_0 \r F(-3)
\]
such that the composition
\[
G_2\oplus F_1 \xrightarrow{\begin{pmatrix} f_2& \gamma_1\\ 0 & -M\end{pmatrix}}
G_1\oplus F_0\r F(-3)
\]
is zero.

Introducing the explicit form of $(G_i)_i$, $(f_i)_i$ given by \eqref{ref-4.2-38},
and using Proposition \ref{ref-3.1.1-19} we find that a map $J \r F(-3)$ is  given by a pair of maps 
\[
\mu:A(-v-1)\r F(-3)
\]
\[
\beta:A(-v-2)^{a_{v+2}}\r F(-3)
\]
(identified with scalar vectors as in \eqref{ref-4.1-37}) such that the composition
\[
A(-v-2)^2\oplus A(-v-3)^{b_{v+3}} \xrightarrow{
\begin{pmatrix}
-l_2 & l_1 & \gamma_1(v+1,v+3)\\
0 & 0 & -M(v+2,v+3)
\end{pmatrix}
}
A(-v-1)\oplus A(-v-2)^{a_{v+2}} \xrightarrow{
\begin{pmatrix} 
\mu & \beta
\end{pmatrix}
} F
\]
is zero.

Let $p$ be the point associated to $F$. Since $(l_1)_p=(l_2)_p=0$ we obtain the conditions
\begin{equation}
\label{ref-4.13-52}
\begin{pmatrix}
\mu & \beta
\end{pmatrix}
\begin{pmatrix}
\gamma_1(v+1,v+3)_p\\M(v+2,v+3)_p
\end{pmatrix}=0
\end{equation}
To use this we have to know what $\gamma_1(v+1,v+3)$ is. From the commutative diagram \eqref{ref-4.3-40} we obtain the identity
\[
\rho \cdot \gamma_1(v+1,v+3)=\lambda\cdot M(u,v+3)
\]
where $\lambda=\gamma_0(u,u)$. Evaluation in $p$ yields
\[
\gamma_1(v+1,v+3)_p=\lambda\cdot M(u,v+3)_p
\]
so that \eqref{ref-4.13-52} is equivalent to 
\[
\begin{pmatrix}
\mu & \beta
\end{pmatrix}
\begin{pmatrix}
\lambda\cdot M(u,v+3)_p\\M(v+2,v+3)_p
\end{pmatrix}=0
\]
Now $
\begin{pmatrix} 
\lambda \cdot M(u,v+3)_p \\ 
M(v+2,v+3)_p
\end{pmatrix}
$ is a $(a_{v+2}+1)\times b_{v+3}$ matrix. Since $b_{v+3} < a_{v+2}+1$ (Proposition \ref{ref-3.1.1-19}) we would expect \eqref{ref-4.13-52} to have $a_{v+2}+1-b_{v+3}$ independent solutions. To have more, $
\begin{pmatrix}
\lambda \cdot M(u,v+3)_p \\
M(v+2,v+3)_p 
\end{pmatrix}
$ has to have non-maximal rank. I.e.\ there should be a non-zero solution to the equation
\[
\begin{pmatrix}
\lambda\cdot M(u,v+3)_p\\M(v+2,v+3)_p
\end{pmatrix}\cdot \delta=0
\]
which may be broken up into two sets of equations
\begin{equation}
\label{ref-4.14-53}
\lambda\cdot M(u,v+3)_p\cdot  \delta=0
\end{equation}
\begin{equation}
\label{ref-4.15-54}
M(v+2,v+3)_p\cdot \delta=0
\end{equation}
and we also still have
\begin{equation}
\label{ref-4.16-55}
\lambda\cdot M(u,u+1)_p=0
\end{equation}
We view \eqref{ref-4.14-53}\eqref{ref-4.15-54} and \eqref{ref-4.16-55} as a system of $1+a_{v+2}+b_{u+1}$ equations in the variety $\PP^{a_{u}-1} \times \PP^2 \times \PP^{b_{v+3}-1}$. Since (Condition C)
\[
1 + a_{v+2} + b_{u+1} = \dim (\PP^{a_{u}-1} \times \PP^2 \times \PP^{b_{v+3}-1}) = a_u + b_{v+3}
\]
the existence of a solution can be decided numerically. The intersection product we have to compute is 
\[
(r+s+t)(s+t)^{a_{v+2}}(r+s)^{b_{u+1}}
\]
This product contains the term
\[
s^2 t^{a_{v+2}-1}r^{b_{u+1}}
\]
which is non-zero in the Chow ring (using Condition C). \qed
\end{case}
\def\qed{}
\end{proof} 

\subsection{Estimating the dimension of { $\Ext^1_{\hat{A}}(\hat{J},\hat{J})$}}
In this section we prove the following result.
\begin{propositions} 
Assume that Condition C holds. Let $I \in H_\varphi$ be generic and let $J$ be as in Condition D. Then
\begin{equation}
\label{ref-4.17-56}
\dim_{k} \Ext^1_{\hat{A}}(\hat{J},\hat{J})\le
\begin{cases}
\dim H_\varphi+a_u-1 & \text{if $v=u$}\\
\dim H_\varphi+a_u+1-b_{u+1} & \text{if $v\ge u+1$}
\end{cases}
\end{equation}
\end{propositions}
\begin{proof}
It has been shown in \cite{DV2} that $H_\varphi$ is the moduli-space of ideals in $A$ of projective dimension one which have Hilbert series $\varphi$. Let $\tilde{I}\subset A_{H_\varphi}$ be the corresponding universal bundle. Let $\Mscr$ be the moduli-space of pairs $(J,I)$ such that $I\in H_\varphi$ and $h_J=\psi$. To show that $\Mscr$ exists on may realize it as a closed subscheme of 
\[
\underline{\Proj}\, S_{H_\varphi}(\tilde{I}_u\oplus\cdots \oplus
\tilde{I}_v)
\]

Sending $(J,I)$ to $I$ defines a map $q:\Mscr\r H_\varphi$. We have
an exact sequence
\begin{equation}
\label{ref-4.18-57}
0\r T_{(J,I)}q^{-1}I\r T_{(J,I)}\Mscr\r T_{I} H_\varphi
\end{equation}
Assume now that $I$ is generic and put $\Vscr=q^{-1}I$ as above. By Proposition \ref{ref-4.3.1-45} we know that $\Vscr$ is smooth. Hence
\[
\dim T_{(J,I)} \Mscr \le  \dim \Vscr + \dim H_\varphi
\]
Applying Proposition \ref{ref-4.3.1-45} again, it follows that for $I$ generic the dimension of $T_{(J,I)}\Mscr$ is bounded by the right hand side of \eqref{ref-4.17-56}.

Since $\Ext^1_{\hat{A}}(\hat{J},\hat{J})$ is the tangent space of $\Mscr$ at $(J,I)$ for $\hat{J}=(J \,\,\, I)$ this finishes the proof.
\end{proof}
\begin{remarks} 
It is not hard to see that \eqref{ref-4.17-56} is actually an equality. This follows from the easily proved fact that
the map $q$ is generically smooth.
\end{remarks}

\subsection{Tying things together}
Combining the results of the previous two sections we see that if Condition C holds we have for a suitable
choice of  $J$
\[
\dim_{k} \Ext^1_A(J,J) - \dim_{k} \Ext^1_{\hat{A}}(\hat{J},\hat{J})
\ge
\begin{cases}
\dim H_\psi-\dim H_\varphi+a_{v+3}-a_u+1&\text{if $v=u$}\\
\dim H_\psi-\dim H_\varphi+a_{v+2}-b_{v+3}-a_u+b_{u+1}&\text{if $v=u+1$}\\
\dim H_\psi-\dim H_\varphi+a_{v+2}-b_{v+3}-a_u+b_{u+1}+1&\text{if $v\ge u+2$}
\end{cases}
\]
We may combine this with Proposition \ref{ref-3.2.1-20} which works out as (using Proposition \ref{ref-3.1.1-19}) 
\[
\dim H_\psi-\dim H_\varphi=
\begin{cases}
a_u+b_{v+3}-1 &\text{if $v=u$}\\
a_u-b_{u+1}-a_{v+2}+b_{v+3}+1&\text{if $v=u+1$}\\
a_u-b_{u+1}-a_{v+2}+b_{v+3}&\text{if $v\ge u+2$}
\end{cases}
\]
We then obtain 
\[
\dim_{k} \Ext^1_A(J,J)-\dim_{k} \Ext^1_{\hat{A}}(\hat{J},\hat{J})
\ge
\begin{cases}
b_{v+3}&\text{if $v=u$} \\
1&\text{if $v\ge u+1$}
\end{cases}
\]
Hence in all cases we obtain a strictly positive result. This finishes the proof that Condition C implies Condition D.
\begin{remarks} 
As in Remark \ref{ref-3.4.1-33} it is possible to prove directly the converse implication $D \Rightarrow C$.
\end{remarks}

\section{The implication $D \Rightarrow A$}
In this section $(\varphi,\psi)$ will have the same meaning as in \S\ref{ref-3-14} and we also keep the associated notations. We assume that Condition D holds. Let $I$ be a graded ideal corresponding to a generic point in $H_\varphi$. According to Condition D there exists an ideal $J\subset I$ with $h_J=\psi$ such that there is an $\eta \in \Ext^1_A(J,J)$ which is not in the image of $\Ext^1_{\hat{A}}(\hat{J},\hat{J})$.

We identify $\eta$ with a one parameter deformation $J'$ of $J$. I.e.\ $J'$ is a flat $A[\epsilon]$-module where $\epsilon^2=0$ such that $J'\otimes_{k[\epsilon]}k\cong J$ and such that the short exact sequence
\[
0\r J \xrightarrow{\epsilon\cdot} J'\r J \r 0
\]
corresponds to $\eta$.

In \S\ref{ref-4.2-39} we have written $J$ as the homology of a complex. It follows for example from (the dual version of) \cite[Thm 3.9]{lowen3}, or directly, that $J'$ is the homology of a complex of the form
\begin{equation}
\label{ref-5.1-58}
0 \r G_3[\epsilon] \xrightarrow{
\begin{pmatrix}
f'_3 \\ 
P\epsilon
\end{pmatrix}
} G_2[\epsilon] \oplus F_1[\epsilon]
\xrightarrow{
\begin{pmatrix}
f'_2 & \gamma'_1 \\
Q\epsilon & -M'
\end{pmatrix}
}
G_1[\epsilon]\oplus F_0[\epsilon]
\xrightarrow{
\begin{pmatrix}
f'_1 & \gamma'_0
\end{pmatrix}}
G_0[\epsilon]\r 0
\end{equation}
where for a matrix $U$ over $A$, $U'$ means a lift of $U$ to $A[\epsilon]$. Recall that $G_3=A(-v-3)$.
\begin{lemma} \label{ref-5.1-59}
We have $P(v+3,v+3)\neq 0$.
\end{lemma}
\begin{proof} 
Assume on the contrary $P(v+3,v+3)=0$. Using Proposition \ref{ref-3.1.1-19} it follows that $P$ has its image in 
$F_{11}=\oplus_{j\le u+1}A(-j)^{b_{j}}$.

The fact that \eqref{ref-5.1-58} is a complex implies that $Qf_3=MP$. Thus we have a commutative diagram
\[
\begin{CD}
0 @>>> G_3 @>{f_{3}}>> G_2 @>{f_2}>> G_1 @>{f_1}>> G_0  \\
@. @VP_1VV  @VVQV\\ 
@. F_{11}@>M_{11}>> F_0\\
@. @VP_2VV @|\\
@. F_1 @>>M> F_0
\end{CD}
\] 
where $P_2$ is the inclusion and $M_{11}=MP_2$, $P=P_2P_1$. Put
\[
D=\coker(F_{11}\r F_1)
\]

Then $(P_1,Q)$ represents an element of $\Ext^2_A(F,D)=\Ext^1_A(D,F(-3))^\ast=0$, where the last equality is for
degree reasons.

It follows that there exist maps
\begin{align*}
R:G_1&\r F_0\\
T_1:G_2&\r F_{11}
\end{align*}
such that 
\begin{align*}
Q&=Rf_2+M_{11}T_1\\
P_1&=T_1f_3
\end{align*}
Putting $T=P_2T_1$ we obtain
\begin{align*}
Q&=Rf_2+MT\\
P&=Tf_3
\end{align*}
We can now construct the following lifting of the commutative diagram
\eqref{ref-4.5-42}:
\[
\xymatrix{
0 \ar[r] & G_3[\epsilon] \ar[r]^(0.4)
{\begin{pmatrix}f'_3\\ P\epsilon\end{pmatrix}}
 & G_2[\epsilon]\oplus F_1[\epsilon] \ar[d]_{\begin{pmatrix}T\epsilon & -1\end{pmatrix}}
\ar[rr]^{
\begin{pmatrix}
f'_2 & \gamma'_1\\
Q\epsilon & -M'
\end{pmatrix}
}
&&
G_1[\epsilon]\oplus F_0[\epsilon]\ar[d]^{\begin{pmatrix} -R\epsilon & 1\end{pmatrix}}
\ar[rr]^(0.6){\begin{pmatrix}f'_1 & \gamma'_0\end{pmatrix}}
&&G_0[\epsilon] \ar[r]& 0\\
&0\ar[r]& F_1[\epsilon] \ar[rr]_{M'+R\gamma_1\epsilon} &&F_0[\epsilon]\ar[r]& 0
}
\]
Taking homology we see that there is a first order deformation $I'$ of $I$ together with a lift of the inclusion $J\r I$ to a map $J'\r I'$. But this contradicts the assumption that $\eta$ is not in the image of $\Ext_{\hat{A}}^1(\hat{J},\hat{J})$.
\end{proof}
In particular, Lemmma \ref{ref-5.1-59} implies that $b_v+3 \neq 0$.
It will now be convenient to rearrange \eqref{ref-5.1-58}. Using the previous lemma and the fact that the rightmost map in \eqref{ref-5.1-58} is split it follows that $J'$ has a free resolution of the form
\[
0\r G_3[\epsilon]\xrightarrow{
\begin{pmatrix}
\epsilon\\
\alpha_0+\alpha_1\epsilon
\end{pmatrix}
}
G_3[\epsilon]\oplus H_1[\epsilon]
\xrightarrow{
\begin{pmatrix}
\beta_0+\beta_1\epsilon & \delta_0+\delta_1\epsilon
\end{pmatrix}
}
H_0[\epsilon]\r J' \r 0
\]
which leads to the following equations
\begin{align*}
\delta_0\alpha_0&=0\\
\beta_0+\delta_1\alpha_0+\delta_0\alpha_1&=0
\end{align*}
Using these equations we can construct the following complex $C_t$ over $A[t]$
\[
0\r G_3[t]\xrightarrow{
\begin{pmatrix}
t\\
\alpha_0+\alpha_1t
\end{pmatrix}
}
G_3[t]\oplus H_1[t]
\xrightarrow{
\begin{pmatrix}
\beta_0-\delta_1\alpha_1 t& \delta_0+\delta_1 t
\end{pmatrix}
}
H_0[t]
\]
For $\theta\in k$ put $C_\theta=C\otimes_{k[t]}k[t]/(t-\theta)$. Clearly $C_0$ is a resolution of $J$. By semi-continuity we find that for all but a finite number of $\theta$, $C_\theta$ is the resolution of a rank one $A$-module $J_\theta$. Furthermore we have $J_0=J$ and $\pd J_\theta=1$ for $\theta\neq 0$. 

\medskip

Let $\Jscr_\theta$ be the rank one $\Oscr_{\PP^2}$-module corresponding to $J_\theta$.  $\Jscr_\theta$ represents a point of $H_\psi$. Since $I/J$ has finite length, $J_0=J$ and $I$ define the same object in $\Hilb_n(\PP^2)$. Hence we have constructed a one parameter family of objects in $\Hilb_n(\PP^2)$ connecting a generic object in $H_\varphi$ to an object in $H_\psi$. This shows that indeed $H_\varphi$ is in the closure of $H_\psi$.

\appendix
\section{Hilbert graphs} \label{ref-A-59} 

For low values of $n$ the Hilbert graph is rather trivial. But when $n$ becomes bigger the number of Hilbert functions increase rapidly (see Remark \ref{ref-1.3-3}) and so the Hilbert graphs become more complicated. As an illustration we have included the Hilbert graph for $n = 17$ where we used Theorem \ref{ref-1.5-8} to solve the incidence problems of length zero (the picture gives no information on more complicated incidence problems). By convention the minimal Hilbert
series is on top. 

The reader will notice that the Hilbert graph contains pentagons. This shows that the Hilbert graph is not catenary and also contradicts \cite[Lemme 2.1.2]{Guerimand}.

\newpage

\unitlength 1mm
\begin{picture}(100.50,200.50)(0,0)

\linethickness{0.15mm}
\put(90.00,180.00){\line(0,1){10.00}}

\linethickness{0.15mm}
\put(90.00,170.00){\line(0,1){10.00}}

\linethickness{0.15mm}
\multiput(90.00,160.00)(0.12,-0.12){83}{\line(1,0){0.12}}

\linethickness{0.15mm}
\multiput(90.00,140.00)(0.12,0.12){83}{\line(1,0){0.12}}

\linethickness{0.15mm}
\multiput(80.00,150.00)(0.12,-0.12){83}{\line(1,0){0.12}}

\linethickness{0.15mm}
\multiput(70.00,140.00)(0.12,-0.12){83}{\line(1,0){0.12}}

\linethickness{0.15mm}
\multiput(80.00,130.00)(0.12,-0.12){83}{\line(1,0){0.12}}

\linethickness{0.15mm}
\multiput(60.00,130.00)(0.12,-0.12){83}{\line(1,0){0.12}}

\linethickness{0.15mm}
\multiput(70.00,120.00)(0.12,-0.12){83}{\line(1,0){0.12}}

\linethickness{0.15mm}
\multiput(80.00,110.00)(0.12,-0.12){83}{\line(1,0){0.12}}

\linethickness{0.15mm}
\multiput(50.00,60.00)(0.12,-0.12){83}{\line(1,0){0.12}}

\linethickness{0.15mm}
\multiput(60.00,50.00)(0.12,-0.12){83}{\line(1,0){0.12}}

\linethickness{0.15mm}
\multiput(40.00,50.00)(0.12,0.12){83}{\line(1,0){0.12}}

\linethickness{0.15mm}
\put(70.00,20.00){\line(0,1){20.00}}

\linethickness{0.15mm}
\multiput(60.00,30.00)(0.12,-0.12){83}{\line(1,0){0.12}}

\linethickness{0.15mm}
\multiput(60.00,10.00)(0.12,0.12){83}{\line(1,0){0.12}}

\linethickness{0.15mm}
\multiput(60.00,10.00)(0.12,-0.12){83}{\line(1,0){0.12}}

\linethickness{0.15mm}
\multiput(70.00,0.00)(0.12,0.12){83}{\line(1,0){0.12}}

\linethickness{0.15mm}
\multiput(90.00,160.00)(0,1.82){6}{\line(0,1){0.91}}

\linethickness{0.15mm}
\multiput(80.00,150.00)(1.33,1.33){8}{\multiput(0,0)(0.11,0.11){6}{\line(1,0){0.11}}}

\linethickness{0.15mm}
\multiput(70.00,140.00)(1.33,1.33){8}{\multiput(0,0)(0.11,0.11){6}{\line(1,0){0.11}}}

\linethickness{0.15mm}
\multiput(80.00,130.00)(1.33,1.33){8}{\multiput(0,0)(0.11,0.11){6}{\line(1,0){0.11}}}

\linethickness{0.15mm}
\multiput(90.00,160.00)(0,1.82){6}{\line(0,1){0.91}}

\linethickness{0.15mm}
\multiput(60.00,130.00)(1.33,1.33){8}{\multiput(0,0)(0.11,0.11){6}{\line(1,0){0.11}}}

\linethickness{0.15mm}
\multiput(70.00,120.00)(1.33,1.33){8}{\multiput(0,0)(0.11,0.11){6}{\line(1,0){0.11}}}

\linethickness{0.15mm}
\multiput(80.00,110.00)(1.33,1.33){8}{\multiput(0,0)(0.11,0.11){6}{\line(1,0){0.11}}}

\linethickness{0.15mm}
\put(90.00,80.00){\line(0,1){20.00}}

\linethickness{0.15mm}
\multiput(50.00,120.00)(0.12,-0.12){333}{\line(1,0){0.12}}

\linethickness{0.15mm}
\multiput(50.00,120.00)(1.33,1.33){8}{\multiput(0,0)(0.11,0.11){6}{\line(1,0){0.11}}}

\linethickness{0.15mm}
\multiput(60.00,110.00)(1.33,1.33){8}{\multiput(0,0)(0.11,0.11){6}{\line(1,0){0.11}}}

\linethickness{0.15mm}
\multiput(70.00,100.00)(1.33,1.33){8}{\multiput(0,0)(0.11,0.11){6}{\line(1,0){0.11}}}

\linethickness{0.15mm}
\multiput(60.00,90.00)(1.33,1.33){8}{\multiput(0,0)(0.11,0.11){6}{\line(1,0){0.11}}}

\linethickness{0.15mm}
\multiput(70.00,80.00)(1.33,1.33){8}{\multiput(0,0)(0.11,0.11){6}{\line(1,0){0.11}}}

\linethickness{0.15mm}
\multiput(40.00,110.00)(0.12,0.12){83}{\line(1,0){0.12}}

\linethickness{0.15mm}
\multiput(50.00,100.00)(0.12,0.12){83}{\line(1,0){0.12}}

\linethickness{0.15mm}
\multiput(40.00,110.00)(0.12,-0.12){83}{\line(1,0){0.12}}

\linethickness{0.15mm}
\multiput(50.00,100.00)(1.33,-1.33){8}{\multiput(0,0)(0.11,-0.11){6}{\line(1,0){0.11}}}

\linethickness{0.15mm}
\multiput(70.00,80.00)(1.33,-1.33){8}{\multiput(0,0)(0.11,-0.11){6}{\line(1,0){0.11}}}

\linethickness{0.15mm}
\multiput(60.00,90.00)(0.12,-0.12){83}{\line(1,0){0.12}}

\linethickness{0.15mm}
\multiput(80.00,70.00)(1.33,1.33){8}{\multiput(0,0)(0.11,0.11){6}{\line(1,0){0.11}}}

\linethickness{0.15mm}
\multiput(50.00,80.00)(1.33,-1.33){8}{\multiput(0,0)(0.11,-0.11){6}{\line(1,0){0.11}}}

\linethickness{0.15mm}
\multiput(60.00,70.00)(1.33,1.33){8}{\multiput(0,0)(0.11,0.11){6}{\line(1,0){0.11}}}

\linethickness{0.15mm}
\multiput(60.00,70.00)(0.12,-0.12){83}{\line(1,0){0.12}}

\linethickness{0.15mm}
\multiput(70.00,60.00)(0.12,0.12){83}{\line(1,0){0.12}}

\linethickness{0.15mm}
\multiput(80.00,70.00)(0.12,-0.12){83}{\line(1,0){0.12}}

\linethickness{0.15mm}
\multiput(50.00,80.00)(0.12,0.12){83}{\line(1,0){0.12}}

\linethickness{0.15mm}
\multiput(50.00,60.00)(1.33,1.33){8}{\multiput(0,0)(0.11,0.11){6}{\line(1,0){0.11}}}

\linethickness{0.15mm}
\multiput(60.00,50.00)(1.33,1.33){8}{\multiput(0,0)(0.11,0.11){6}{\line(1,0){0.11}}}

\linethickness{0.15mm}
\multiput(60.00,50.00)(1.94,0.65){16}{\multiput(0,0)(0.32,0.11){3}{\line(1,0){0.32}}}

\linethickness{0.15mm}
\multiput(40.00,50.00)(1.33,-1.33){8}{\multiput(0,0)(0.11,-0.11){6}{\line(1,0){0.11}}}

\linethickness{0.15mm}
\multiput(50.00,40.00)(1.33,1.33){8}{\multiput(0,0)(0.11,0.11){6}{\line(1,0){0.11}}}

\linethickness{0.15mm}
\multiput(50.00,40.00)(1.33,-1.33){8}{\multiput(0,0)(0.11,-0.11){6}{\line(1,0){0.11}}}

\linethickness{0.15mm}
\multiput(70.00,20.00)(1.33,-1.33){8}{\multiput(0,0)(0.11,-0.11){6}{\line(1,0){0.11}}}

\put(20.63,185.00){\makebox(0,0)[cc]{means $H_{\varphi} \subset \overline{H_{\psi}}$}}

\linethickness{0.15mm}
\put(0.00,180.00){\line(0,1){10.00}}

\linethickness{0.15mm}
\multiput(0.00,165.00)(0,1.82){6}{\line(0,1){0.91}}

\linethickness{0.15mm}
\put(0.00,150.00){\line(0,1){10.00}}

\put(20.63,170.00){\makebox(0,0)[cc]{means $H_{\varphi} \not\subset \overline{H_{\psi}}$}}

\put(28.63,155.00){\makebox(0,0)[cc]{means $(\varphi,\psi)$ has type zero}}

\put(3.00,190.00){\makebox(0,0)[cc]{$\varphi$}}

\put(2.50,180.00){\makebox(0,0)[cc]{$\psi$}}

\put(2.00,155.00){\makebox(0,0)[cc]{$0$}}

\put(53.75,106.88){\makebox(0,0)[cc]{$0$}}

\put(67.50,86.88){\makebox(0,0)[cc]{$0$}}

\put(73.13,66.25){\makebox(0,0)[cc]{$0$}}

\put(72.00,30.00){\makebox(0,0)[cc]{$0$}}

\put(3.13,175.00){\makebox(0,0)[cc]{$\varphi$}}

\put(3.00,160.00){\makebox(0,0)[cc]{$\varphi$}}

\put(3.00,150.00){\makebox(0,0)[cc]{$\psi$}}

\put(3.00,165.00){\makebox(0,0)[cc]{$\psi$}}

\linethickness{0.15mm}
\put(90.00,190.00){\circle*{1.00}}

\linethickness{0.15mm}
\put(90.00,180.00){\circle*{1.00}}

\linethickness{0.15mm}
\put(90.00,170.00){\circle*{1.00}}

\linethickness{0.15mm}
\put(90.00,160.00){\circle*{1.00}}

\linethickness{0.15mm}
\put(80.00,150.00){\circle*{1.00}}

\linethickness{0.15mm}
\put(100.00,150.00){\circle*{1.00}}

\linethickness{0.15mm}
\put(70.00,140.00){\circle*{1.00}}

\linethickness{0.15mm}
\put(90.00,140.00){\circle*{1.00}}

\linethickness{0.15mm}
\put(60.00,130.00){\circle*{1.00}}

\linethickness{0.15mm}
\put(80.00,130.00){\circle*{1.00}}

\linethickness{0.15mm}
\put(50.00,120.00){\circle*{1.00}}

\linethickness{0.15mm}
\put(70.00,120.00){\circle*{1.00}}

\linethickness{0.15mm}
\put(90.00,120.00){\circle*{1.00}}

\linethickness{0.15mm}
\put(40.00,110.00){\circle*{1.00}}

\linethickness{0.15mm}
\put(60.00,110.00){\circle*{1.00}}

\linethickness{0.15mm}
\put(80.00,110.00){\circle*{1.00}}

\linethickness{0.15mm}
\put(50.00,100.00){\circle*{1.00}}

\linethickness{0.15mm}
\put(70.00,100.00){\circle*{1.00}}

\linethickness{0.15mm}
\put(90.00,100.00){\circle*{1.00}}

\linethickness{0.15mm}
\put(60.00,90.00){\circle*{1.00}}

\linethickness{0.15mm}
\put(80.00,90.00){\circle*{1.00}}

\linethickness{0.15mm}
\put(50.00,80.00){\circle*{1.00}}

\linethickness{0.15mm}
\put(70.00,80.00){\circle*{1.00}}

\linethickness{0.15mm}
\put(90.00,80.00){\circle*{1.00}}

\linethickness{0.15mm}
\put(60.00,70.00){\circle*{1.00}}

\linethickness{0.15mm}
\put(80.00,70.00){\circle*{1.00}}

\linethickness{0.15mm}
\put(50.00,60.00){\circle*{1.00}}

\linethickness{0.15mm}
\put(70.00,60.00){\circle*{1.00}}

\linethickness{0.15mm}
\put(90.00,60.00){\circle*{1.00}}

\linethickness{0.15mm}
\put(40.00,50.00){\circle*{1.00}}

\linethickness{0.15mm}
\put(60.00,50.00){\circle*{1.00}}

\linethickness{0.15mm}
\put(50.00,40.00){\circle*{1.00}}

\linethickness{0.15mm}
\put(70.00,40.00){\circle*{1.00}}

\linethickness{0.15mm}
\put(60.00,30.00){\circle*{1.00}}

\linethickness{0.15mm}
\put(60.00,30.00){\circle*{1.00}}

\linethickness{0.15mm}
\put(70.00,20.00){\circle*{1.00}}

\linethickness{0.15mm}
\put(60.00,10.00){\circle*{1.00}}

\linethickness{0.15mm}
\put(80.00,10.00){\circle*{1.00}}

\linethickness{0.15mm}
\put(70.00,0.00){\circle*{1.00}}

\linethickness{0.15mm}
\put(0.00,190.00){\circle*{1.00}}

\linethickness{0.15mm}
\put(0.00,180.00){\circle*{1.00}}

\linethickness{0.15mm}
\put(0.00,175.00){\circle*{1.00}}

\linethickness{0.15mm}
\put(0.00,165.00){\circle*{1.00}}

\linethickness{0.15mm}
\put(0.00,160.00){\circle*{1.00}}

\linethickness{0.15mm}
\put(0.00,150.00){\circle*{1.00}}

\end{picture}

\section{A visual criterion for incidence problems of length zero}
\label{ref-B-60}
In this appendix we provide a visual criterion for the three conditions in Theorem \ref{ref-1.5-8}. The reader may easily check these using Condition C, Proposition \ref{ref-3.1.1-19} and \eqref{ref-3.1-15}. As before we let $(\varphi, \psi)$ be a pair of Hilbert series of degree $n$ and length zero. Then $H_{\varphi} \subset \overline{H_{\psi}}$ if and only if the Castelnuovo diagram $s_{\varphi}$ of $\varphi$ has one of the following forms, where the diagram $s_{\psi}$ is obtained by moving the upper square as indicated.

\unitlength 0.8mm
\begin{picture}(190.00,188.13)(15,0)

\linethickness{0.15mm}
\put(25.00,100.00){\line(1,0){5.00}}

\linethickness{0.15mm}
\multiput(30.00,90.00)(0,1.82){6}{\line(0,1){0.91}}

\linethickness{0.15mm}
\put(40.00,90.00){\line(1,0){5.00}}

\linethickness{0.15mm}
\put(45.00,85.00){\line(0,1){5.00}}

\linethickness{0.15mm}
\put(40.00,85.00){\line(1,0){5.00}}

\linethickness{0.15mm}
\multiput(45.00,75.00)(0,1.82){6}{\line(0,1){0.91}}

\linethickness{0.15mm}
\put(45.00,75.00){\line(1,0){5.00}}

\linethickness{0.15mm}
\put(20.00,90.00){\line(0,1){5.00}}

\linethickness{0.15mm}
\put(15.00,90.00){\line(1,0){5.00}}

\linethickness{0.15mm}
\put(15.00,85.00){\line(0,1){5.00}}

\linethickness{0.15mm}
\put(10.00,85.00){\line(1,0){5.00}}

\linethickness{0.15mm}
\put(10.00,80.00){\line(0,1){5.00}}

\linethickness{0.15mm}
\put(25.00,95.00){\line(0,1){5.00}}

\linethickness{0.15mm}
\put(20.00,95.00){\line(1,0){5.00}}

\linethickness{0.15mm}
\multiput(5.00,75.00)(1.43,1.43){4}{\multiput(0,0)(0.12,0.12){6}{\line(1,0){0.12}}}

\linethickness{0.15mm}
\put(30.00,90.00){\line(1,0){10.00}}

\linethickness{0.15mm}
\put(40.00,85.00){\line(0,1){5.00}}

\linethickness{0.15mm}
\qbezier(42.50,92.50)(41.88,108.75)(34.38,95.00)

\linethickness{0.15mm}
\put(34.38,95.00){\line(0,1){1.88}}

\linethickness{0.15mm}
\multiput(34.38,95.00)(0.19,0.13){10}{\line(1,0){0.19}}

\linethickness{0.15mm}
\put(50.00,75.00){\line(1,0){5.00}}

\linethickness{0.15mm}
\multiput(55.00,75.00)(1.82,0){6}{\line(1,0){0.91}}

\linethickness{0.15mm}
\put(80.00,90.00){\line(0,1){10.00}}
\put(80.00,100.00){\vector(0,1){0.12}}
\put(80.00,90.00){\vector(0,-1){0.12}}

\linethickness{0.15mm}
\put(45.00,65.00){\line(1,0){20.00}}
\put(65.00,65.00){\vector(1,0){0.12}}
\put(45.00,65.00){\vector(-1,0){0.12}}

\put(55.00,60.00){\makebox(0,0)[cc]{$\geq 2$}}

\put(85.00,95.00){\makebox(0,0)[cc]{$\geq 0$}}

\put(85.00,82.50){\makebox(0,0)[cc]{$\geq 1$}}

\linethickness{0.15mm}
\put(30.00,65.00){\line(1,0){15.00}}
\put(45.00,65.00){\vector(1,0){0.12}}
\put(30.00,65.00){\vector(-1,0){0.12}}

\put(37.50,60.00){\makebox(0,0)[cc]{$3$}}

\linethickness{0.15mm}
\put(80.00,75.00){\line(0,1){15.00}}
\put(80.00,90.00){\vector(0,1){0.12}}
\put(80.00,75.00){\vector(0,-1){0.12}}

\linethickness{0.15mm}
\put(130.00,110.00){\line(1,0){5.00}}

\linethickness{0.15mm}
\multiput(135.00,100.00)(0,1.82){6}{\line(0,1){0.91}}

\linethickness{0.15mm}
\put(145.00,100.00){\line(1,0){5.00}}

\linethickness{0.15mm}
\put(150.00,95.00){\line(0,1){5.00}}

\linethickness{0.15mm}
\put(145.00,95.00){\line(1,0){5.00}}

\linethickness{0.15mm}
\put(125.00,100.00){\line(0,1){5.00}}

\linethickness{0.15mm}
\put(120.00,100.00){\line(1,0){5.00}}

\linethickness{0.15mm}
\put(120.00,95.00){\line(0,1){5.00}}

\linethickness{0.15mm}
\put(115.00,95.00){\line(1,0){5.00}}

\linethickness{0.15mm}
\put(115.00,90.00){\line(0,1){5.00}}

\linethickness{0.15mm}
\put(130.00,105.00){\line(0,1){5.00}}

\linethickness{0.15mm}
\put(125.00,105.00){\line(1,0){5.00}}

\linethickness{0.15mm}
\multiput(100.00,75.00)(1.43,1.43){4}{\multiput(0,0)(0.12,0.12){6}{\line(1,0){0.12}}}

\linethickness{0.15mm}
\put(135.00,100.00){\line(1,0){10.00}}

\linethickness{0.15mm}
\put(145.00,95.00){\line(0,1){5.00}}

\linethickness{0.15mm}
\put(139.38,105.00){\line(0,1){1.88}}

\linethickness{0.15mm}
\multiput(139.38,105.00)(0.19,0.13){10}{\line(1,0){0.19}}

\linethickness{0.15mm}
\put(180.00,100.00){\line(0,1){10.00}}
\put(180.00,110.00){\vector(0,1){0.12}}
\put(180.00,100.00){\vector(0,-1){0.12}}

\put(185.00,105.00){\makebox(0,0)[cc]{$\geq 0$}}

\linethickness{0.15mm}
\put(135.00,65.00){\line(1,0){15.00}}
\put(150.00,65.00){\vector(1,0){0.12}}
\put(135.00,65.00){\vector(-1,0){0.12}}

\put(142.50,60.00){\makebox(0,0)[cc]{$3$}}

\linethickness{0.15mm}
\multiput(150.00,85.00)(0,1.82){6}{\line(0,1){0.91}}

\linethickness{0.15mm}
\put(150.00,80.00){\line(0,1){5.00}}

\linethickness{0.15mm}
\put(150.00,80.00){\line(1,0){5.00}}

\linethickness{0.15mm}
\put(155.00,75.00){\line(0,1){5.00}}

\linethickness{0.15mm}
\put(155.00,75.00){\line(1,0){5.00}}

\linethickness{0.15mm}
\multiput(160.00,75.00)(1.82,0){6}{\line(1,0){0.91}}

\linethickness{0.15mm}
\put(180.00,80.00){\line(0,1){20.00}}
\put(180.00,100.00){\vector(0,1){0.12}}
\put(180.00,80.00){\vector(0,-1){0.12}}

\linethickness{0.15mm}
\put(155.00,65.00){\line(1,0){15.00}}
\put(170.00,65.00){\vector(1,0){0.12}}
\put(155.00,65.00){\vector(-1,0){0.12}}

\put(162.50,60.00){\makebox(0,0)[cc]{$\geq 1$}}

\put(160.00,60.00){\makebox(0,0)[cc]{}}

\put(185.00,90.00){\makebox(0,0)[cc]{$\geq 2$}}

\put(182.50,90.00){\makebox(0,0)[cc]{}}

\linethickness{0.15mm}
\put(110.00,90.00){\line(1,0){5.00}}

\linethickness{0.15mm}
\put(110.00,85.00){\line(0,1){5.00}}

\linethickness{0.15mm}
\put(105.00,85.00){\line(1,0){5.00}}

\linethickness{0.15mm}
\put(105.00,80.00){\line(0,1){5.00}}

\linethickness{0.15mm}
\put(115.00,45.00){\line(1,0){5.00}}

\linethickness{0.15mm}
\multiput(120.00,35.00)(0,1.82){6}{\line(0,1){0.91}}

\linethickness{0.15mm}
\put(145.00,35.00){\line(1,0){5.00}}

\linethickness{0.15mm}
\put(150.00,30.00){\line(0,1){5.00}}

\linethickness{0.15mm}
\put(145.00,30.00){\line(1,0){5.00}}

\linethickness{0.15mm}
\multiput(150.00,20.00)(0,1.82){6}{\line(0,1){0.91}}

\linethickness{0.15mm}
\put(150.00,20.00){\line(1,0){5.00}}

\linethickness{0.15mm}
\put(110.00,35.00){\line(0,1){5.00}}

\linethickness{0.15mm}
\put(105.00,35.00){\line(1,0){5.00}}

\linethickness{0.15mm}
\put(105.00,30.00){\line(0,1){5.00}}

\linethickness{0.15mm}
\put(100.00,30.00){\line(1,0){5.00}}

\linethickness{0.15mm}
\put(100.00,25.00){\line(0,1){5.00}}

\linethickness{0.15mm}
\put(115.00,40.00){\line(0,1){5.00}}

\linethickness{0.15mm}
\put(110.00,40.00){\line(1,0){5.00}}

\linethickness{0.15mm}
\multiput(95.00,20.00)(1.43,1.43){4}{\multiput(0,0)(0.12,0.12){6}{\line(1,0){0.12}}}

\linethickness{0.15mm}
\put(145.00,30.00){\line(0,1){5.00}}

\linethickness{0.15mm}
\put(155.00,20.00){\line(1,0){5.00}}

\linethickness{0.15mm}
\multiput(160.00,20.00)(1.82,0){6}{\line(1,0){0.91}}

\linethickness{0.15mm}
\put(180.00,35.00){\line(0,1){10.00}}
\put(180.00,45.00){\vector(0,1){0.12}}
\put(180.00,35.00){\vector(0,-1){0.12}}

\linethickness{0.15mm}
\put(150.00,10.00){\line(1,0){20.00}}
\put(170.00,10.00){\vector(1,0){0.12}}
\put(150.00,10.00){\vector(-1,0){0.12}}

\put(160.00,5.00){\makebox(0,0)[cc]{$\geq 2$}}

\put(185.00,40.00){\makebox(0,0)[cc]{$\geq 0$}}

\put(185.00,27.50){\makebox(0,0)[cc]{$\geq 1$}}

\linethickness{0.15mm}
\put(180.00,20.00){\line(0,1){15.00}}
\put(180.00,35.00){\vector(0,1){0.12}}
\put(180.00,20.00){\vector(0,-1){0.12}}

\linethickness{0.15mm}
\put(140.00,35.00){\line(1,0){5.00}}

\linethickness{0.15mm}
\put(120.00,35.00){\line(1,0){10.00}}

\linethickness{0.15mm}
\multiput(130.00,35.00)(1.82,0){6}{\line(1,0){0.91}}

\linethickness{0.15mm}
\put(120.00,10.00){\line(1,0){30.00}}
\put(150.00,10.00){\vector(1,0){0.12}}
\put(120.00,10.00){\vector(-1,0){0.12}}

\put(135.00,5.00){\makebox(0,0)[cc]{$\geq 4$}}

\linethickness{0.15mm}
\qbezier(147.50,37.50)(134.38,58.13)(125.00,40.00)

\linethickness{0.15mm}
\put(125.00,40.00){\line(0,1){1.25}}

\linethickness{0.15mm}
\multiput(125.00,40.00)(0.25,0.13){5}{\line(1,0){0.25}}

\linethickness{0.15mm}
\multiput(20.00,40.00)(0,1.82){6}{\line(0,1){0.91}}

\linethickness{0.15mm}
\put(30.00,35.00){\line(1,0){5.00}}

\linethickness{0.15mm}
\put(35.00,30.00){\line(0,1){5.00}}

\linethickness{0.15mm}
\put(30.00,30.00){\line(1,0){5.00}}

\linethickness{0.15mm}
\multiput(35.00,20.00)(0,1.82){6}{\line(0,1){0.91}}

\linethickness{0.15mm}
\put(35.00,20.00){\line(1,0){5.00}}

\linethickness{0.15mm}
\put(20.00,35.00){\line(1,0){10.00}}

\linethickness{0.15mm}
\put(30.00,30.00){\line(0,1){5.00}}

\linethickness{0.15mm}
\qbezier(32.50,37.50)(31.88,53.75)(24.38,40.00)

\linethickness{0.15mm}
\put(24.38,40.00){\line(0,1){1.88}}

\linethickness{0.15mm}
\multiput(24.38,40.00)(0.19,0.13){10}{\line(1,0){0.19}}

\linethickness{0.15mm}
\put(40.00,20.00){\line(1,0){5.00}}

\linethickness{0.15mm}
\multiput(45.00,20.00)(1.82,0){6}{\line(1,0){0.91}}

\linethickness{0.15mm}
\put(35.00,10.00){\line(1,0){20.00}}
\put(55.00,10.00){\vector(1,0){0.12}}
\put(35.00,10.00){\vector(-1,0){0.12}}

\put(45.00,5.00){\makebox(0,0)[cc]{$\geq 2$}}

\put(70.00,27.50){\makebox(0,0)[cc]{$\geq 1$}}

\linethickness{0.15mm}
\put(20.00,10.00){\line(1,0){15.00}}
\put(35.00,10.00){\vector(1,0){0.12}}
\put(20.00,10.00){\vector(-1,0){0.12}}

\put(27.50,5.00){\makebox(0,0)[cc]{$3$}}

\linethickness{0.15mm}
\put(65.00,20.00){\line(0,1){15.00}}
\put(65.00,35.00){\vector(0,1){0.12}}
\put(65.00,20.00){\vector(0,-1){0.12}}

\linethickness{0.15mm}
\put(20.00,35.00){\line(0,1){5.00}}

\linethickness{0.15mm}
\put(15.00,50.00){\line(1,0){5.00}}

\linethickness{0.15mm}
\put(10.00,50.00){\line(1,0){5.00}}

\linethickness{0.15mm}
\multiput(0.00,50.00)(1.82,0){6}{\line(1,0){0.91}}

\linethickness{0.15mm}
\put(65.00,35.00){\line(0,1){15.00}}
\put(65.00,50.00){\vector(0,1){0.12}}
\put(65.00,35.00){\vector(0,-1){0.12}}

\put(70.00,43.75){\makebox(0,0)[cc]{$\geq 1$}}

\put(67.50,43.75){\makebox(0,0)[cc]{}}

\linethickness{0.15mm}
\put(0.00,10.00){\line(1,0){20.00}}
\put(20.00,10.00){\vector(1,0){0.12}}
\put(0.00,10.00){\vector(-1,0){0.12}}

\put(10.00,5.00){\makebox(0,0)[cc]{$\geq 2$}}

\linethickness{0.15mm}
\qbezier(147.50,102.50)(146.88,118.75)(139.38,105.00)

\linethickness{0.15mm}
\put(35.00,170.00){\line(1,0){5.00}}

\linethickness{0.15mm}
\put(40.00,170.00){\line(1,0){5.00}}

\linethickness{0.15mm}
\multiput(45.00,160.00)(0,1.82){6}{\line(0,1){0.91}}

\linethickness{0.15mm}
\put(45.00,160.00){\line(1,0){5.00}}

\linethickness{0.15mm}
\put(50.00,155.00){\line(0,1){5.00}}

\linethickness{0.15mm}
\put(45.00,155.00){\line(0,1){5.00}}

\linethickness{0.15mm}
\put(45.00,155.00){\line(1,0){5.00}}

\linethickness{0.15mm}
\multiput(50.00,145.00)(0,1.82){6}{\line(0,1){0.91}}

\linethickness{0.15mm}
\put(50.00,145.00){\line(1,0){5.00}}

\linethickness{0.15mm}
\multiput(55.00,135.00)(0,1.82){6}{\line(0,1){0.91}}

\linethickness{0.15mm}
\put(55.00,135.00){\line(1,0){5.00}}

\linethickness{0.15mm}
\multiput(60.00,135.00)(1.82,0){6}{\line(1,0){0.91}}

\linethickness{0.15mm}
\put(30.00,160.00){\line(0,1){5.00}}

\linethickness{0.15mm}
\put(25.00,160.00){\line(1,0){5.00}}

\linethickness{0.15mm}
\put(25.00,155.00){\line(0,1){5.00}}

\linethickness{0.15mm}
\put(20.00,155.00){\line(1,0){5.00}}

\linethickness{0.15mm}
\put(20.00,150.00){\line(0,1){5.00}}

\linethickness{0.15mm}
\put(35.00,165.00){\line(0,1){5.00}}

\linethickness{0.15mm}
\put(30.00,165.00){\line(1,0){5.00}}

\linethickness{0.15mm}
\put(15.00,150.00){\line(1,0){5.00}}

\linethickness{0.15mm}
\put(15.00,145.00){\line(0,1){5.00}}

\linethickness{0.15mm}
\put(10.00,145.00){\line(1,0){5.00}}

\linethickness{0.15mm}
\put(10.00,140.00){\line(0,1){5.00}}

\linethickness{0.15mm}
\multiput(5.00,135.00)(1.43,1.43){4}{\multiput(0,0)(0.12,0.12){6}{\line(1,0){0.12}}}

\linethickness{0.15mm}
\put(80.00,160.00){\line(0,1){10.00}}
\put(80.00,170.00){\vector(0,1){0.12}}
\put(80.00,160.00){\vector(0,-1){0.12}}

\put(85.00,165.00){\makebox(0,0)[cc]{$\geq 0$}}

\linethickness{0.15mm}
\put(80.00,145.00){\line(0,1){15.00}}
\put(80.00,160.00){\vector(0,1){0.12}}
\put(80.00,145.00){\vector(0,-1){0.12}}

\linethickness{0.15mm}
\put(80.00,135.00){\line(0,1){10.00}}
\put(80.00,145.00){\vector(0,1){0.12}}
\put(80.00,135.00){\vector(0,-1){0.12}}

\linethickness{0.15mm}
\put(55.00,125.00){\line(1,0){15.00}}
\put(70.00,125.00){\vector(1,0){0.12}}
\put(55.00,125.00){\vector(-1,0){0.12}}

\put(87.50,152.50){\makebox(0,0)[cc]{$C \geq 1$}}

\put(87.50,140.00){\makebox(0,0)[cc]{$D \geq 0$}}

\put(99.38,146.25){\makebox(0,0)[cc]{where $C > D$}}

\put(62.50,120.00){\makebox(0,0)[cc]{$\geq 1$}}

\linethickness{0.15mm}
\qbezier(50.63,162.50)(56.25,188.13)(43.75,175.00)

\linethickness{0.15mm}
\multiput(43.75,175.00)(0.13,0.38){5}{\line(0,1){0.38}}

\linethickness{0.15mm}
\multiput(43.75,175.00)(0.38,0.13){5}{\line(1,0){0.38}}

\linethickness{0.15mm}
\put(35.00,125.00){\line(1,0){10.00}}
\put(45.00,125.00){\vector(1,0){0.12}}
\put(35.00,125.00){\vector(-1,0){0.12}}

\put(40.00,120.00){\makebox(0,0)[cc]{$2$}}

\linethickness{0.15mm}
\put(130.00,153.75){\line(1,0){5.00}}

\linethickness{0.15mm}
\multiput(135.00,143.75)(0,1.82){6}{\line(0,1){0.91}}

\linethickness{0.15mm}
\put(135.00,143.75){\line(1,0){5.00}}

\linethickness{0.15mm}
\put(140.00,138.75){\line(0,1){5.00}}

\linethickness{0.15mm}
\put(135.00,138.75){\line(0,1){5.00}}

\linethickness{0.15mm}
\put(135.00,138.75){\line(1,0){5.00}}

\linethickness{0.15mm}
\multiput(140.00,128.75)(0,1.82){6}{\line(0,1){0.91}}

\linethickness{0.15mm}
\put(140.00,128.75){\line(1,0){5.00}}

\linethickness{0.15mm}
\multiput(145.00,118.75)(0,1.82){6}{\line(0,1){0.91}}

\linethickness{0.15mm}
\put(145.00,120.00){\line(1,0){5.00}}

\linethickness{0.15mm}
\multiput(150.00,120.00)(1.82,0){6}{\line(1,0){0.91}}

\linethickness{0.15mm}
\put(170.00,143.75){\line(0,1){10.00}}
\put(170.00,153.75){\vector(0,1){0.12}}
\put(170.00,143.75){\vector(0,-1){0.12}}

\put(175.00,148.75){\makebox(0,0)[cc]{$\geq 0$}}

\linethickness{0.15mm}
\put(170.00,128.75){\line(0,1){15.00}}
\put(170.00,143.75){\vector(0,1){0.12}}
\put(170.00,128.75){\vector(0,-1){0.12}}

\linethickness{0.15mm}
\put(170.00,118.75){\line(0,1){10.00}}
\put(170.00,128.75){\vector(0,1){0.12}}
\put(170.00,118.75){\vector(0,-1){0.12}}

\put(177.50,136.25){\makebox(0,0)[cc]{$C \geq 1$}}

\put(177.50,123.75){\makebox(0,0)[cc]{$D \geq 0$}}

\put(190.00,130.00){\makebox(0,0)[cc]{where $C > D$}}

\linethickness{0.15mm}
\qbezier(140.63,146.25)(146.25,171.88)(133.75,158.75)

\linethickness{0.15mm}
\multiput(133.75,158.75)(0.13,0.38){5}{\line(0,1){0.38}}

\linethickness{0.15mm}
\multiput(133.75,158.75)(0.38,0.13){5}{\line(1,0){0.38}}

\linethickness{0.15mm}
\put(130.00,153.75){\line(0,1){5.00}}

\linethickness{0.15mm}
\multiput(130.00,158.75)(0,1.82){6}{\line(0,1){0.91}}

\linethickness{0.15mm}
\put(125.00,168.75){\line(1,0){5.00}}

\linethickness{0.15mm}
\multiput(125.00,168.75)(0,1.82){6}{\line(0,1){0.91}}

\linethickness{0.15mm}
\put(120.00,178.75){\line(1,0){5.00}}

\linethickness{0.15mm}
\put(170.00,168.75){\line(0,1){10.00}}
\put(170.00,178.75){\vector(0,1){0.12}}
\put(170.00,168.75){\vector(0,-1){0.12}}

\linethickness{0.15mm}
\put(170.00,153.75){\line(0,1){15.00}}
\put(170.00,168.75){\vector(0,1){0.12}}
\put(170.00,153.75){\vector(0,-1){0.12}}

\put(177.50,173.75){\makebox(0,0)[cc]{$A \geq 0$}}

\put(177.50,161.25){\makebox(0,0)[cc]{$B \geq 1$}}

\put(190.00,167.50){\makebox(0,0)[cc]{where $A < B$}}

\linethickness{0.15mm}
\multiput(110.00,178.75)(1.82,0){6}{\line(1,0){0.91}}

\end{picture}

\newpage


\ifx\undefined\bysame
\newcommand{\bysame}{\leavevmode\hbox to3em{\hrulefill}\,}
\fi

\end{document}